\documentclass[12pt]{amsart}
\usepackage{amssymb,amscd}
\input{"xiamacro.sty"}

\begin{document}
\begin{abstract}
We place the representation variety in the broader context of abelian and nonabelian cohomology.  We outline  the equivalent constructions of the moduli spaces of  flat bundles, of smooth integrable connections,  and of holomorphic integrable connections over a compact K\"ahler manifold.  In addition, we describe  the moduli space of Higgs bundles and how it relates to the representation variety.
We attempt to avoid abstraction, but strive to
present and clarify the unifying ideas underlying the theory.
\end{abstract}

\title{Abelian and non-abelian cohomology}

\author[Xia]{Eugene Z. Xia}
\address{
National Center for Theoretical Sciences\\
Department of Mathematics\\
National Cheng-kung University\\
Tainan 701, Taiwan \\
\tt{ezxia@ncku.edu.tw}
({\it Xia})}

\date{\today}

\thanks{The author gratefully acknowledges partial support by the National Science Council, Taiwan with grants 97-2115-M-006-001-MY3 and 100-2115-M-006-004-MY2 and the Institute of Mathematical Sciences of the National University of Singapore.  He also thanks the latter's hospitality during his visit.}

\maketitle
\tableofcontents

\section{Prelude}\label{sec:Prelude}
At first glance, the representation variety is a strange object.  From the perspective of group theory, the fundamental group of a surface does not seem to have anything to do with Lie groups such as $\SL(2,\R)$.  Yet there is a deep connection from uniformization:  the group $\SL(2,\R)$ acts on the symmetric space $\SL(2,\R)/\U(1)$ which is identified with the upper half of the complex plane $\C$.  Given a discrete subgroup $\Gamma < \SL(2,\R)$, the quotient $X = \Gamma\backslash\SL(2,\R)/\U(1)$ is a Riemann surface with fundamental group isomorphic to $\Gamma$.  Thus a discrete representation $\rho$ of the fundamental group of $X$ gives rise to a Riemann (hyperbolic) surface.  Hence representation varieties are related to Teichm\"uler spaces.  This is the main perspective and focus of the summer school.  

There is another and related perspective of the representation variety that may also be traced back to Riemann and which is just as important.  Roughly speaking, one is interested in the holomorphic linear ordinary differential equation
$$\frac{\d^2 f}{\d z^2} + a(z) \frac{\d f}{\d z} + b(z) f = 0,$$
where $a$ and $b$ are rational functions and together with three distinct poles on the Riemann sphere $\P^1(\C)$ under some additional conditions.    
A particular example is the Euler hypergeometric equation
$$
z(1-z)\frac{d^2f}{dz^2} + (c - (a + b + 1)z) \frac{df}{dz} - abf = 0
$$
with $a, b, c \in \C$.
Then the solutions are the classical Gauss hypergeometric functions.  These are complicated functions with multiple indices that can make any non-analyst dizzy and which satisfy still more complicated Kummer relations.

However under mild conditions and by a standard elementary transformation, one may rewrite this second order equation as a system of first order equations
\begin{equation}\label{eq:Riemann-Hilbert}
\frac{\d F}{\d z} + \frac{A_0 F}{z} + \frac{A_1 F}{z-1} = 0, \ \ A_i \in \gl(2,\C), \ i = 0,1,
\end{equation}
where $\gl(2,\C)$ consists of $2 \times 2$ complex matrices and is the Lie algebra of $\GL(2,\C)$.  Let $A_\infty = -A_0 -A_1$.  If $A_0 = A_1 = 0$, then $F$ is a meromorphic function on $\P^1(\C)$ with poles at $\{0, 1, \infty\}$ such that $\frac{\d F}{\d z} = 0$.  This implies that $F$ is a constant function. In general, such a solution $F$ only exists locally.  Globally, $F$ exists as a multivalued function with monodromy at $\{0, 1, \infty\}$.  This gives rise to the concept of a local system (system of locally constant functions) and a monodromy representation
$$\rho : \pi_1(\P^1(\C) \setminus \{0, 1, \infty\}) \lto \GL(2,\C).$$  Riemann pointed out (in our modern language) that the $\GL(2,\C)$-representation variety of a 3-holed sphere is a point.  In the language of hyperbolic geometry, the hyperbolic structure on a pair of pants is uniquely determined by the lengths of its three geodesic boundary components when $\im(\rho) < \SL(2,\R)$.  In other words, the solutions to the equation are determined by the eigenvalues of $A_1, A_2$ and $A_\infty$ only.  These eigenvalues are easy to compute and together they determine the monodromy matrices $C_0, C_1$ and $C_\infty$ uniquely up to simultaneous conjugation equivalence, i.e. $C_i$ is equal to a conjugate of $\exp(A_i)$ and $C_0 C_1 C_\infty = \id$.  This set of $C_i$'s determines $\rho$.  Riemann then concluded that the hypergeometric functions corresponding to the solutions of the equation are uniquely determined by the corresponding monodromy representation $\rho$.  In this way, he read off the Kummer relations among the hypergeometric functions without much calculation \cite{Ka1, Ri1}.  The representation varieties into other Lie groups have similar interpretations in terms of holomorphic differential equations.


It is then of great interest to obtain as much information as possible concerning the representation varieties.
To understand a space means to discover natural and canonical topological and geometric structures on that space.  Let $X$ be a manifold and $G$ a Lie group and $\MB$ the resulting representation variety (the Betti moduli space) of $X$.  A priori, $\MB$ is a set, one point for each representation up to conjugation equivalence.  However the algebraic/analytic structure on $G$ gives $\MB$ a variety structure, making it a representation variety.  If $G$ is reductive and $X$ is compact and symplectic, then the Poincar\'e duality and a nondegenerate symmetric bi-invariant form on $G$ give a symplectic structure on $\MB$, making it a symplectic manifold or orbifold.

One method to discover topological and geometric structures is to give alternative constructions that give rise to new topological and geometric structures.  A representation corresponds to a flat $G$-vector bundle on $X$.  A flat vector bundle is one where (locally) constant sections make sense.  In order for this to happen, the transition maps over the overlapping covering open sets of $X$ must be constant.  These constant transition maps then give rise to a representation.  Hence the representation variety corresponds to the moduli space of flat $G$-vector bundles.  This gives rise to the \v{C}ech moduli space $\MC$ and is the point of view of local systems.

A function $f$ is (locally) constant if and only if 
its partial derivatives vanish everywhere.  Differentiation 
is $\R$-linear
and satisfies the Leibniz rule when applied to products of functions.
Differentiation
is an example of more general first order linear differential operators, namely, vector fields that also satisfy $\R$-linearity and the Leibniz rule. This allows one to say that $f$ is {\em constant} if and only if $Vf = 0$ for each vector field $V$.    The set of vector fields $\T$ forms an (infinite dimensional) vector space over $\R$.  The very fact that vector fields operate on functions implies that $\T$, considered as an (additive) abelian group, acts on the space of functions, making the latter a $\T$-module.
One then considers the family of $\T$-modules where the $\T$-action satisfies the $\R$-linearity condition and the Leibniz rule.  This leads to the concept of {\em connection}.  In such a module, one may then ask whether there is a reasonable notion of constant $e$ defined as $V e = 0$ for all $V \in \T$.
This further leads to the idea of {\em integrability} and $\D$-modules.
On a smooth manifold $X$, vector fields are well-defined;
hence, the notion of integrable connections may be extended to smooth manifolds.  Hence the \v{C}ech moduli space has an equivalent construction as the moduli space of sheaves of $\D$-modules or moduli space of flat $G$-connections.  This is the smooth de Rham moduli space $\MdR$.

Suppose $X$ is a complex manifold.
Then a function $f$ on $X$ is constant if it is holomorphic and all its holomorphic partial derivatives vanish.  This gives rise to the holomorphic de Rham moduli space, parameterizing holomorpic structures with integrable holomorphic connections.  This is the perspective of holomorphic differential equations (See (\ref{eq:Riemann-Hilbert}) above).  The holomorphic connections are further generalized to holomorphic $\lambda$-connections for $\lambda \in \C$.  This is an algebraic/complex analytic construction and the resulting moduli space $\MdR^\lambda$ has a natural variety structure.  If $\lambda \in \C^\times$, then all the moduli spaces paramaterized by $\lambda$ have isomorphic underlying variety structure.  When $X$ is compact, $\MdR^0$ acquires a {\em different} variety structure while maintaining the same underlying topological/differentiable structure as that of $\MdR^1$.  This $\MdR^0$ is the Dolbeault moduli space.  This new variety structure is compatible with that on the Betti moduli space.  These two variety structures together with the symplectic structure give rise to a hyperk\"ahler structure on the underlying space.  To see applications of this structure in the study of geometric structures on $X$, see the survey paper \cite{G0}.

Section \ref{sec:Hom} constructs the representation variety, the Betti moduli space $\MB$, for a manifold and gives two simple examples and points
out the subtle issues of the construction.  Section \ref{sec:groupoid} is a quick and informal introduction to groupoids.  It first serves an accounting purpose to keep track of the various objects we are dealing with.  More importantly, with the possible exception of the Betti moduli space, the constructions of the other moduli spaces as varieties are technically daunting.  Rather, this set of introductory notes only describes the groupoids of the objects these moduli spaces parameterize and the relations between the objects among these groupoids.

Sections~\ref{sec:local R-functions} and \ref{sec:local C-functions} deal with local constructions of the smooth and holomorphic de Rham spaces on the unit ball.
Section~\ref{sec: local examples} gives some examples of connections and especially important is the example that shows that a connection must be integrable for it to define constant functions.

Section \ref{sec: smooth manifolds} is a brief outline of the concepts of smooth and complex manifolds in preparation for global constructions on manifolds.  Section \ref{sec: sheaves} describes the smooth, flat and holomorphic bundles on a manifold.  Sections \ref{sec: global smooth} and \ref{sec:global holomorphic} run in parallel to Sections \ref{sec:local R-functions} and \ref{sec:local C-functions} and construct the flat connections, holomorphic structures and flat holomorphic $\lambda$-connections
respectively.

Section \ref{sec:iso} describes the relations between the \v{C}ech, Betti and de Rham moduli spaces and how they relate to the Dolbeault moduli space in the case of compact K\"ahler manifolds.  As we are leaving out the actual constructions of these moduli spaces, we state the equivalence and give examples of these equivalences in the rank-1 trivial cases in Section \ref{sec:simple examples}.  The final Section \ref{sec:final remarks} hints on the complex structure arising from the Dolbeault moduli space and how this complex structure and the one from the Betti construction and the symplectic structure give rise to a hyperk\"ahler structure on the smooth part of the representation variety.

The subject is structure rich.
These notes attempt to cover only the most basic ideas--the leitmotif is the simple observation that a smooth function is (locally) constant if and only if all its partial derivatives vanish.  They are also designed to complement and supplement \cite{GX1} which focuses on the rank-1 sheaves on Riemann surfaces.   These notes attempt to expand the scope to high rank sheaves over K\"ahler manifolds, but focus more on detailed local descriptions.  While \cite{GX1} is succinct and rigorous, these notes are informal and avoid abstraction.  The emphasis is on presenting the basic unifying ideas underlying the theory.



\section{The representation variety}\label{sec:Hom}
Let $X$ be a manifold and $G$ the general linear group $\GL(r,\C)$.
Let $x \in X$ and $\pi = \pi_1(X, x)$ be the topological fundamental group of $X$.  We assume that $\pi$ is finitely generated with $N$ generators and $M$ relations.
Then the space of representations $\Hom(\pi, G)$ is a subvariety of $G^N$ defined by the $M$ relations.  The group $G$ acts on $\Hom(\pi, G)$ by conjugation:
$$
G \times \Hom(\pi, G) \lto \Hom(\pi, G), \ \ (g, \rho) \mapsto g \rho g^{-1}.
$$
One would like to form a quotient of this $G$-action.  However
the orbits of the $G$-conjugation action may not be closed in $\Hom(\pi, G)$.
This implies that the natural geometric quotient (orbit space) of the $G$-conjugation action may not be Hausdorff.
To resolve this and obtain a Hausdorff quotient, one must either remove the non-closed orbits or identify them with their closure.

Since $G$ is the rank-$r$ linear group,
each representation $\rho \in \Hom(\pi, G)$ induces a rank-$r$ representation of $\pi$ by composition.  We say $\rho$ is reductive if every $\pi$-invariant subspace has a $\pi$-invariant complement in $\C^{\oplus r}$.  Denote by $\Hom(\pi, G)^+$ the subspace of reductive representations.  The $G$-conjugation action preserves $\Hom(\pi, G)^+$ and with closed orbits.  Denote by $\MB$ the quotient $\Hom(\pi, G)^+/G$.  Then $\MB$ parameterizes the equivalence classes of reductive representations.

Alternatively, $\Hom(\pi, G)$ is an affine subvariety of $G^N$, defined by the $M$ equations with its resulting coordinate ring $R$.
Hence the $G$-action on $\Hom(\pi_1, G)$ induces a $G$-action on $R$.  Denote by $R^G$ the invariant subring.  Then $\MB$ is the affine variety defined by the ring $R^G$.

\noindent {\em Example:}  Let $X$ be the compact orientable surface of genus $g$ and $G = \GL(1,\C) = \C^\times$.
Set $[A,B] = ABA^{-1}B^{-1}.$  Then
$$\pi = \langle A_i, B_i, 1 \le i \le g \ |  \ \prod_{i=1}^g [A_i,B_i] \rangle$$
Since $G$ is abelian,
$\Hom(\pi_1, \C^\times) \cong  (\C^\times)^{2g}$ and the $G$-action is trivial.  Hence
$\MB \cong (\C^\times)^{2g}$.

\noindent {\em Example:} Let $X$ be the compact orientable surface of genus $g$ and $G = \GL(2,\C)$. Then
$$\Hom(\pi_1, G) \cong
\{a_i,b_i \in G : 1 \le i \le g, \prod_{i=1}^g [a_i, b_i] = e \in G\}.$$
Consider the representation $\rho$ corresponding to
$$
a_i, b_i = \left[
\begin{array}{cc}
1 & 1 \\
0 & 1
\end{array}
\right], \ \ 1 \le i \le g.
$$
The $G$-orbit of $\rho$ is not closed because it does not contain the trivial representation.  Hence the corresponding geometric quotient space is not Hausdorff.

\section{A reluctant tour of category and groupoid}\label{sec:groupoid}
The basic problem is to classify a family of objects up to isomorphism and impose a natural and universal geometry over the isomorphism classes.   This problem is often cast in the language of category and groupoid.  We will take a somewhat informal approach and refer
to (Section 1, \cite{GX1}) for a careful introduction of these concepts.  For us,
\begin{defin}
A {\bf groupoid} $\CC$ is a category consisting of
\begin{enumerate}

\item A set of objects $\Obj(\CC)$.

\item A set (quite often a group) $\G$ of equivalence relations on the objects.
\end{enumerate}
\end{defin}

We often use $\CC$ to denote both the groupoid and its set of objects $\Obj(\CC)$ as well.  The equivalence relation is often described by a group action $\G \times \CC \lto \CC$.
From the groupoid, one obtains a set $\Iso(\CC)$ of isomorphism classes and then attempts to construct a natural geometric structure on $\Iso(\CC)$ in a universal way.  This is the concept of moduli space.

In our specific construction,
\begin{defin}
The {\bf Betti groupoid} $\Hom(\pi, G)$ consists of
\begin{enumerate}

\item The set of objects: $\Hom(\pi, G)$.

\item The equivalence group $G$-action:
$$
G \times \Hom(\pi, G) \lto \Hom(\pi, G), \ \ (g, \rho) \mapsto g\rho g^{-1}.
$$
\end{enumerate}
\end{defin}
From these two ingredients, one obtains the orbit family $\Iso(\Hom(\pi,G))$.  The set $\Hom(\pi, G)$
acquires a variety structure from $G$ and $\MB$ parameterizes the reductive representations in $\Iso(\Hom(\pi,G))$.

The Betti groupoid $\Hom(\pi,G)$ is our concern.  To understand this space better, we present a few other related groupoids, namely the \v{C}ech, de Rham and Dolbeault groupoids.  See Definitions ~\ref{def:Cech groupoid smooth}, ~\ref{def:Cech groupoid}, ~\ref{def: smooth de Rham groupoid} and \ref{def:holomorphic de Rham groupoid} for their constructions.

\section{Local functions and their derivatives over $\R^n$}\label{sec:local R-functions}
\subsection{Functions}
Let $X$ be the unit ball in $\R^n$ with coordinates $x = (x_1, \cdots, x_n)$ and suppose that $f : X \lto \C$ is a smooth (infinitely differentiable) function.  A simple observation is that $f$ is a constant function if and only if
$$\frac{\partial f}{\partial x_i} := \frac{\partial \Re(f)}{\partial x_i} + {\bf i} \frac{\partial \Im(f)}{\partial x_i}= 0$$ for all $ 1\le i \le n$; here ${\bf i} \in \C$ with ${\bf i}^2 = -1$.  This simple observation is the kernel of the theory.
Denote by $\Os$ the set of smooth $\C$-valued functions on $X$ and by $\Os_c$ the subset of constant functions.  These two sets form commutative rings with $\Os_c \cong \C$ and $\Os_c \hookrightarrow \Os$ a ring homomorphism.

\subsection{Exterior differentiation}
The linear first order differential operator $\frac{\partial}{\partial x_i}$ is an example of vector fields on $X$.  Denote by $\T$ the space of all such vector fields.  If $f \in \Os$ is a function and $V \in \T$ is a vector field, then $fV$ is again a vector field.  This means that the space $\T$ is a left $\Os$-module.    A function $f \in \Os$ is constant if and only if $V(f) = 0$ for each vector field $V \in \T$.  This motivates the definition of the exterior differentiation.  Let $\Omega^1$ be the space of 1-forms (covector fields) on $X$.  For each smooth function $f : X \to \C$, $d f$ is a 1-form:
$$
d : \Os \lto \Omega^1  \ \ :\ \   df(V) = V(f), \ \ \forall V \in \T,
$$
i.e. $df$ takes a vector field $V$ and returns a function:  $V(f)$.  Then $f \in \Os_c$ if and only if $df(V) = 0$ for each $V \in \T$, or equivalently, $df = 0 \in \Omega^1$ as a covector field.  In other words, the ring $\Os_c$ may be defined as
\begin{equation}
\Os^d = \{f \in \Os : df = 0\}.
\end{equation}
From this, one constructs the $\Os$-modules of $i$-forms $\Omega^i = \wedge_{j=1}^i \Omega^1$ and generalizes $d$ to obtain the standard de Rham complex
$$
(\Omega^\bullet, d) :  \Omega^0 \stackrel{d}{\lto} \Omega^{1} \cdots \Omega^i \stackrel{d}{\lto} \Omega^{i+1} \cdots,
$$
such that $d^2 = 0$ and satisfying the (generalized) Leibniz rule
$$d(a \wedge b) = da \wedge b + (-1)^i a \wedge db, \ \ \ a \in \Omega^i.$$

\subsection{Connections}
Let
$
\Em = \{ f : X \lto \C^{\oplus r}\}.
$
Then $\Em$ is a free $\Os$-module isomorphic to $\Os^{\oplus r}$.
Let $\Em_c \subseteq \Em$ be the $\Os_c$-module isomorphic to $\Os_c^{\oplus r} \subseteq \Os^{\oplus r}$.
The exterior differentiation is now generalized to an operator
$$D : \Em \lto \Omega^1(\Em) := \Omega^1 \otimes \Em, \ \ D(f) = D(\oplus_{i=1}^r f_i) = \oplus_{i=1}^r d f_i.$$
Now we further generalize the notion of exterior differentiation by enforcing linearity and the Leibniz rule.
\begin{defin}
Let $\Omega^i(\Em) := \Omega^i \otimes \Em$.
An $\Os_c$-linear operator
$$\nabla : \Em \lto \Omega^1(\Em)$$
is a {\bf connection} if
$\nabla(f e) = df \otimes e + f \nabla(e)$ for $f \in \Os$ and $e \in \Em$.
Denote by $\CC_d$ the set of all connections on $\Em$.

For each $i$, a connection $\nabla$ extends to an $\Os_c$-linear map
$$
\nabla : \Omega^i(\Em) \lto \Omega^{i+1}(\Em)
$$
satisfying the generalized Leibniz rule
$$\nabla(\eta \otimes e) = d \eta \otimes e + (-1)^i \eta \wedge \nabla(e).$$
\end{defin}
\begin{rem}
A connection $\nb$ together with a vector field $V \in \T$ induce a map $$\nb_V : \Em \stackrel{\nb}{\lto} \Omega^1(\Em) \stackrel{V}{\lto} \Em.$$  Hence $\nabla$ turns $\Em$ into a $\T$-module, satisfying $\Os_c$-linearity and the Leibniz rule.
\end{rem}
Notice that $D \in \CC_d$.
Suppose $\nabla$ and $\nabla_0$ are two connections.  Then for $f \in \Os, e \in \Em$,
$$
(\nb - \nb_0)(f e) = df \otimes e + f \nb e - df \otimes e - f \nb_0 e = f (\nb - \nb_0)(e).
$$
Hence
$\eta = \nabla - \nabla_0$ is $\Os$-linear and a connection is of the form $\nabla = D + A$,
where
$$A \in \Omega^1(\End(\Em)) := \Omega^1 \otimes \End(\Em)$$
and the set of connections identifies with $\Omega^1(\End(\Em))$.
One may think of $A$ as an $r \times r$ matrix with 1-forms for entries.
Then for $f \in \Os, e \in \Em$,
$$\nabla^2 (f e) = \nabla (df \otimes e + f \nabla (e)) = d^2f \otimes e - df \wedge \nabla(e) + df \wedge \nabla(e) + f\nabla^2(e) = f\nabla^2(e).$$
This implies that the operator
$\nabla^2$
is $\Os$-linear, hence, $\nabla^2 \in \Omega^2(\End(\Em))$.
Again one may think of elements in $\Omega^2(\End(\Em))$ as matrices with scalar 2-form entries.
\begin{defin}
The {\bf curvature operator} is a map
$$
\Fsf : \CC_d \lto \Omega^2(\End(\Em)), \ \ \Fsf(\nb) = \nb^2.
$$
For $\nb \in \CC_d$, $\Fsf(\nb)$ is the {\bf curvature} of $\nb$.
A connection $\nabla$ is {\bf integrable} if $\Fsf(\nabla) = 0$.
If $\nabla$ is integrable, then $\Em^\nabla = \{e \in \Em : \nabla e = 0 \}$ is the $\Os_c$-module of constant functions with respect to $\nabla$.
Denote by $\FF_d\subseteq \CC_d$ the subspace of integrable connections.
\end{defin}

\begin{rem}\label{rem:E^nb = E_c}
Notice that $D \in \FF_d$ and $\Em^D = \Em_c$.  If we consider a connection $\nb$ to be a generalized $D$ and use $\nb$ to define $\Os_c$-modules of constant functions, then it is not enough that $\nb$ is $\Os_c$-linear and satisfies the Leibniz rule.  $\nb$ must satisfy the integrability condition in addition.
\end{rem}
Identifying $\CC_d$ with $\Omega^1(\End(\Em))$, there is an exact sequence of maps
$$
\FF_d\hookrightarrow \Omega^1(\End(\Em)) \stackrel{\Fsf}{\lto} \Omega^2(\End(\Em)).
$$
Notice that $\FF_d$ is not necessarily a vector space because $\nabla = D + A_1 + A_2$ is not necessarily integrable even if $\nabla_1 = D + A_1$ and $\nabla_2 = D + A_2$ are.  Another way of saying this is that $\Fsf$ is not a linear operator in general.

For an integrable $\nabla \in \FF_d$, the $\nb$-de Rham complex is
$$
(\Omega^{\bullet}(\Em), \nb) :  \Omega^{0}(\Em) \stackrel{\nb}{\lto} \Omega^{1}(\Em) \cdots \Omega^{i}(\Em) \stackrel{\nb}{\lto} \Omega^{i+1}(\Em) \cdots.
$$

\begin{rem}
Vector fields in $\T$ are first order differential operators on $\Os$.  By composing these operators, $\T$ generates $\D$, the ring of differential operators.  Integrability of $\nb$ is precisely the needed condition for $\Em$ to inherit a $\D$-module structure from its $\T$-module structure \cite{Bo1, Co2, Sc1}.
\end{rem}

\subsection{The gauge group and its action}
The group $G$ acts on $\C^{\oplus r}$ by matrix-vector multiplication.  Define the $G$-gauge group
$$\G = \Omega^0(G) := \Os(G) = \{g : X \lto G\}.$$
Then $\G$ acts on $\Em$ by $\Os$-module automorphism
$$
\G \times \Em \lto \Em, \ \ (g,e) \mapsto ge, \ \ ge(x) = g(x).e(x),
$$
where $g(x).e(x)$ means matrix-vector multiplication over a point $x \in X$.
Suppose $g \in \G$ and $\nb \in \CC_d$.  Then $g : \Em \stackrel{\cong}{\lto} \Em$ and $\nb$ pulls back to another connection $g^*\nb$ on $\Em$ as follows:
$$
g^*\nb : \Em \stackrel{g}{\lto} \Em \stackrel{\nb}{\lto} \Omega^1(\Em) \stackrel{g^{-1}}{\lto} \Omega^1(\Em).
$$
This gives a $\G$-action on $\CC_d$.

\begin{rem}
The $\Os$-module $\End(\Em)$ should be thought of as the Lie algebra valued functions $$\Os(\g) = \Omega^0(\g) = \{v : X \lto \g\}$$ with the bracket as the commutator, where $\g$ is the Lie algebra of $G$.
\end{rem}
Identify $\CC_d$ with $\Omega^1(\End(\Em)) \cong \Omega^1(\g)$ and let $\nb = D + A$ for $A \in \Omega^1(\End(\Em))$ and $g \in \G$.  Then
$$
g^{-1} \nabla g = D + g^{-1} (Dg) + g^{-1} A g,
$$
$$\Fsf(\nb) = \nb^2 = D A + A \wedge A, \ \ \ \Fsf(g^{-1} \nb g) = g^{-1} \nb^2 g,$$
where $Dg$ and $DA$ mean entry-wise exterior differentiation and $A \wedge A$ means matrix multiplication with entry-wise wedge product.
Hence $\nb^2 = 0$ if and only if $g^{-1}\nb^2g = 0$, i.e. the gauge group action preserves and restricts to the integrable connections
$\G \times \FF_d \lto \FF_d.$
Notice that $\g$ may not be commutative; hence, $A \wedge A$ (more canonically written as $\frac{1}{2}[A,A]$ in the language of differential graded Lie algebra \cite{GM1}) is not necessarily zero.

\begin{defin}
The orbit space $\MdR = \FF_d/\G$ is called the {\bf de Rham moduli space} of $X$.
\end{defin}
The connection $D$ defines a map
$$
\Dsf : \G \lto \CC_d, \ \ \Dsf(g) = g^{-1} D g.
$$
Identifying $\CC_d$ with $\Omega^1(\End(\Em))$, the following sequence of maps is exact
$$
\Omega^0(G) \stackrel{\Dsf}{\lto} \Omega^1(\End(\Em)) \stackrel{\Fsf}{\lto} \Omega^2(\End(\Em)).
$$
Hence the space $\MdR$ is a point.
This is the same as saying that on the unit ball, all integrable connections are $\G$-gauge equivalent.  The general principle is that the curvature is the only local obstruction to a connection being trivial, i.e. isomorphic to $D$. 
\begin{rem}\label{rem: local consant subset smooth}
An integrable connection $\nabla$ produces an $\Os_c$-module $\Em^\nb \cong \Em_c \subseteq \Em$.  The reverse operation is $\Em = \Os \otimes \Em^\nb.$
\end{rem}

\section{Local functions and their derivatives over $\C^n$}\label{sec:local C-functions}
Let $X$ be the unit ball in $\C^n$ with coordinates $z = (z_1, \cdots, z_n)$ and suppose that $f : X \lto \C$ is a smooth function.
\subsection{Holomorphic functions and forms}
As $X$ has a complex structure, one may reformulate the leitmotif as: $f$ is a constant function if and only if $\frac{\partial f}{\partial \bz_i} = 0$ and $\frac{\partial f}{\partial z_i} = 0$ for all $1 \le i \le n$.  In other words, one may define intermediate objects, halfway between the very flexible smooth functions and the very rigid constant functions:
\begin{defin}
A function $f \in \Os$ is {\bf holomorphic} if it satisfies the Cauchy-Riemann equation $\frac{\partial f}{\partial \bar{z_i}} = 0$ for all $1 \le i \le n$.
Denote by $\Os_h$ the ring of holomorphic functions on $X$.
\end{defin}
Then $f$ is constant if and only if $f$ is holomorphic and, in addition, all its holomorphic derivatives are zero, i.e. a constant function $f$ is a holomorphic function that satisfies the additional equations $\frac{\partial f}{\partial z_i} = 0$ for all $i$.  Again notice that $\Os_c \hookrightarrow \Os_h \hookrightarrow \Os$ are ring homomorphisms.

The complex structure on the unit ball $X$ in $\C^n$ gives a decomposition $\T = \T^{1,0} \oplus \T^{0,1}$ as $\Os$-modules, where $\T^{1,0}$ is generated by $\{\frac{\partial}{\partial z_i}: 1 \le i \le n\}$ and $\T^{0,1}$ is generated by $\{\frac{\partial}{\partial \bar{z_i}}: 1 \le i \le n\}$.  This induces a decomposition on the dual space $\Omega^1 = \Omega^{1,0} \oplus \Omega^{0,1}$.  More generally, by taking exterior products of $\Omega^{1,0}$ and $\Omega^{0,1}$, one obtains decompositions and projections
$$\Omega^m = \bigoplus_{q+p=m}\Omega^{q,p}, \ \ \ P_{q,p} : \Omega^m \lto \Omega^{q,p}$$
and the exterior differential operator $d : \Omega^m \lto \Omega^{m+1}$ decomposes as
$$
d = \d + \dbar \ \text{ with } \ \d = P_{q+1,p} \circ d \ \text{ and } \ \ \ \db = P_{q,p+1} \circ d.
$$
Moreover, for fixed $p$ and $q$, the de Rham complex decomposes as
$$
(\Omega^{q,\bullet}, \dbar) :  \Omega^{q,0} \stackrel{\dbar}{\lto} \Omega^{q,1} \cdots \Omega^{q,i} \stackrel{\dbar}{\lto} \Omega^{q,i+1} \cdots
$$
$$
(\Omega^{\bullet,p}, \d) :  \Omega^{0,p} \stackrel{\d}{\lto} \Omega^{1,p} \cdots \Omega^{i,p} \stackrel{\d}{\lto} \Omega^{i+1,p} \cdots,
$$
and
\begin{equation}
\Os_h := \Os^\dbar = \{f \in \Os : \dbar f = 0\}.
\end{equation}
\subsection{The holomorphic structures}
The operator $D$ also decomposes according to types $ D = D' + D'':$
$$D'' : \Em \lto \Omega^{0,1}(\Em) = \Omega^{0,1} \otimes \Em, \ \ D''(f) = D(\oplus_{i=1}^r f_i) = \oplus_{i=1}^r \dbar f_i,$$
$$D' : \Em \lto \Omega^{1,0}(\Em) = \Omega^{1,0} \otimes \Em, \ \ D'(f) = D'(\oplus_{i=1}^r f_i) = \oplus_{i=1}^r \d f_i.$$

\begin{defin}
Let $\Omega^{q,p}(\Em) := \Omega^{q,p} \otimes \Em$.
An almost holomorphic structure on $\Em$ is an operator
$$\bn : \Em \lto \Omega^{0,1}(\Em), \ \ \text{ with } \bn(f e) = \dbar f \otimes e + f \bn(e)$$
for $f \in \Os$ and $e \in \Em$.  Denote by $\Ch$ the space of almost holomorphic structures on $\Em$.

For $p+q=i$, $\bn$ extends to an $\Os_c$-linear map
$$
\bn : \Omega^{p,q}(\Em) \lto \Omega^{p,q+1}(\Em)
$$
satisfying the generalized Leibniz rule
$$\bn(\eta \otimes e) = \dbar \eta \otimes e + (-1)^i \eta \wedge \bn(e).$$
\end{defin}
Define map
$$
\Fsf : \Ch \lto \Omega^{0,2}(\End(\Em)), \ \ \Fsf(\bn) = \bn^2.
$$
\begin{defin}
A {\bf holomorphic structure} $\bn$ on $\Em$ is an almost holomorphic structure satisfying the integrability condition $\Fsf(\bn) = 0$.
Denote by $\Fh$ the set of holomorphic structures on $\Em$.
\end{defin}

Suppose that $\bn$ and $\bn_0$ are two almost holomorphic structures.  Then
$$
(\bn - \bn_0)(f e) = \dbar f \otimes e + f \bn e - \dbar f \otimes e - f \bn_0 e = f (\bn - \bn_0)(e).
$$
Hence $\bn - \bn_0 \in \Omega^{0,1}(\End(\Em))$ and an almost holomorphic structure is of the form $\bn = D'' + A$,
where $A \in \Omega^{0,1}(\End(\Em))$.   Hence the space of almost holomorphic structures $\Ch$ identifies with $\Omega^{0,1}(\End(\Em))$.

When the integrability condition is satisfied, define the Dolbeault complex
$$
(\Omega^{q,\bullet}(\Em), \bn) :  \Omega^{q,0}(\Em) \stackrel{\bn}{\lto} \Omega^{q,1}(\Em) \cdots \Omega^{q,i}(\Em) \stackrel{\bn}{\lto} \Omega^{q,i+1}(\Em) \cdots
$$

Suppose $\bn$ is a holomorphic structure.
Then the $\Os_h$-module $$\Em^\bn = \{e : \bn e = 0\}$$ is the module of holomorphic sections with respect to $\bn$.
Notice that $D'' \in \Fh$.

\subsection{The gauge group action I}
The gauge group $\G$ acts on the almost holomorphic structures, inducing an action
$$
\G \times \Ch \lto \Ch, \ \ (g, \bn) \mapsto g^{-1} \bn g.
$$
This action preserves $\Fh$, hence, restricts to an action
$
\G \times \Fh \lto \Fh.
$
The operator $D''$ also defines a map
$$
\Dsf : \G \lto \Ch, \ \ \Dsf(g) = g^{-1}D''g.
$$
Identifying $\Ch$ with $\Omega^{0,1}(\End(\Em))$, there is an exact sequence of maps
$$
\G \stackrel{\Dsf}{\lto}  \Omega^{0,1}(\End(\Em)) \stackrel{\Fsf}{\lto} \Omega^{0,2}(\End(\Em)).
$$
Hence the space of holomorphic structures $\M_h := \Fh/\G$ is a point.
\begin{rem}\label{rem: local holomorphic subset smooth}
A holomorphic structure $\bn$ produces an $\Os_h$-module $\Em^\bn \cong \Os_h^{\oplus r} \subseteq \Em$.  The reverse operation is $\Em = \Os \otimes \Em^\bn.$
\end{rem}


\subsection{Holomorphic $\lambda$-connections}
At this point, we move to the holomorphic universe (or category) and make a change of notations to emphasize a certain suggestive analogy.  By Remark~\ref{rem: local holomorphic subset smooth}, fix $\Em_h = \Os_h^{\oplus r}$ and let
$$\T_h = \{V \in \T : V(f) \in \Os_h,  \forall f \in \Os_h\}, \ \ \Omh^q = \{\eta \in \Omega^{q,0} : \dbar \eta = 0\},$$  and define the holomorphic de Rham complex
$$
(\Omh^{\bullet}, \d) :  \Omh^{0} \stackrel{\d}{\lto} \Omh^{1} \cdots \Omh^{i} \stackrel{\d}{\lto} \Omh^{i+1} \cdots.
$$
We make a slight generalization:
\begin{defin}
Let
$\Omh^{i}(\Em_h) := \Omh^i \otimes \Em_h$ and
$\lambda \in \C$.
A {\bf holomorphic $\lambda$-connection} or {\bf $\lambda$-connection} on $\Em_h$ is an operator $$\nl : \Em_h \lto \Omh^{1}(\Em_h)$$
satisfying $\nl(f e) = \lambda \d f \otimes e + f \nl(e)$ for $f \in \Os_h$ and $e \in \Em_h$.  Denote by $\Chl$ the space of all $\lambda$-connections on $\Em_h$.

For each $i$, $\nl$ extends to an $\Os_c$-linear map
$$
\nl : \Omh^i(\Em) \lto \Omh^{i+1}(\Em)
$$
satisfying the generalized Leibniz rule
$$\nl(\eta \otimes e) = \lambda \d \eta \otimes e + (-1)^i \eta \wedge \nl(e).$$
\end{defin}
When $\lambda = 1$ we recover the usual definition of connection.
Let
$$
\Fsf : \Chl \lto \Omega_h^2(\End(\Em)), \ \ \Fsf(\nl) = (\nl)^2.
$$
\begin{defin}
A holomorphic $\lambda$-connection $\nl \in \Chl$ is {\bf integrable} if $\Fsf(\nl) = 0$.
Denote by $\Fhl$ the space of integrable $\lambda$-connections on $\Em_h$.
\end{defin}
Suppose that $\nl$ and $\nl_0$ are two $\lambda$-connections on $\Em_h$.  Then
$$
(\nl - \nl_0)(f e) = \lambda \d f \otimes e + f \nl e - \lambda \d f \otimes e - f \nl_0 e = f (\nl - \nl_0)(e).
$$
Hence $\nl - \nl_0 \in \Omega_h^1(\End(\Em_h))$ and a $\lambda$-connection is of the form $\nl = \lambda D' + A$,
where $A \in \Omega_h^1(\End(\Em_h))$.   Hence the space $\Chl$ of $\lambda$-connections identifies with $\Omega_h^1(\End(\Em_h))$.

When the integrability condition is satisfied, the $\lambda$-holomorphic de Rham complex is
$$
(\Omh^{\bullet}(\Em_h), \nl) :  \Omh^{0}(\Em_h) \stackrel{\nl}{\lto} \Omh^{1}(\Em_h) \cdots \Omh^{i}(\Em_h) \stackrel{\nl}{\lto} \Omh^{i+1}(\Em_h) \cdots
$$

If $\nl$ is integrable, then the $\Os_c$-module $\Em_h^{\nl} = \{e : \nl e = 0\}$ is called the module of constant sections with respect to $\nl$.
\begin{rem}
Vector fields in $\T_h$ are first order differential operators on $\Os_h$.
By operator composition, $\T_h$ generates $\D_h$, the ring of differential operators.  Integrability of $\nb^1$ is precisely the needed condition for $\Em_h$ to inherit a $\D_h$-module structure from its $\T_h$-module structure.
\end{rem}

\begin{rem}\label{rem: local lambda connection}
If $\lambda \in \C^\times$, then $\nb^\lambda = \lambda D' + A$ is equivalent to $\nb^1 = D' + \frac{A}{\lambda}$ in the sense that $\nb^\lambda e = 0$ if and only if $\nb^1 e = 0$.

\end{rem}

\subsection{The gauge group action II}
The group $G$ is a complex manifold; hence, it makes sense to say whether a map $g : X \lto G$ is holomorphic.
Notice that each $g \in \G$ is an $\Os$-module automorphism $g : \Em \lto \Em$.
Let
$$
\G_h = \{g \in \G : g : \Em_h \to \Em_h\},
$$
i.e. elements in $\G_h$ restrict to $\Os_h$-module automorphisms of $\Em_h$.
By definition, there is an action
$\G_h \times \Em_h \lto \Em_h$
which induces an action
$$
\Gh \times \Chl \lto \Chl, \ \ (g, \nbl) \mapsto g^{-1} \nbl g = \lambda D' + \lambda g^{-1} D' g + g^{-1} A g,
$$
where $\nbl = \lambda D' + A$.
The $\Gh$-action preserves $\Fhl$.
The operator $D'$ defines a map
$$
\Dsf : \G_h \lto \Chl, \ \  \Dsf(g) = g^{-1} \lambda D' g
$$
and there is an exact sequence of maps
$$
\G_h \stackrel{\Dsf}{\lto} \Chl \stackrel{\Fsf}{\lto} \Omega_h^2(\End(\Em_h)).
$$
The $\lambda$-connections in the same orbits are considered equivalent.
On the unit ball $X$, there is only one integrable $\lambda$-connection on $\Em_h$ up to $\Gh$-gauge equivalence.  Hence $\Fhl/\G_h$
is a point.

\begin{rem}\label{rem: local consant subset holomorphic subset smooth}
The $\Os$-module $\Em$ together with $\bn$ make $\Em^\bn \cong \Em_h$ an $\Os_h$-module.  Furthermore, for $\lambda \in \C^\times$,  $\Em_h$ together with $\nl$ give an $\Os_c$-module isomorphism $\Em_h^\lambda \cong \Em_c$.
To summarize, for $\lambda \in \C^\times$,
$$
\Em_h^\lambda \subseteq \Em_h \subseteq \Em.
$$
In the reverse direction,
$$
\Em_h \cong \Os_h \otimes \Em_c, \ \ \Em \cong \Os \otimes \Em_h, \ \ \Em \cong \Os \otimes \Em_c.
$$
\end{rem}

\section{Local examples in low dimensions}\label{sec: local examples}
\subsection{Connections on an interval} Let $X = (-1, 1) \subset \R$ and suppose that $f : X \lto \C$ is a smooth function.  Then $f$ is a constant function if and only if $\frac{\partial f}{\partial x} = 0$.

Consider the situation of $r = \rk(\Em) = 1$.  Then $G = \C^\times$, $\Em \cong \Os$, $D = d$ and $\CC_d$ identifies with $\Omega^1(\End(\Em)) \cong \Omega^1$.  Let $\nb \in \CC_d$ be a connection.  Then $\nb = D + A = d + A$ and $f \in \Em^\nb$ if and only if
$$\nabla f = df + Af = 0.$$
Choosing $f(0) = C$ and integrating, the solution is
$$f(x) =  Ce^{-\int_0^x A}.$$
The gauge group is
$$\G := \Omega^0(\C^\times) := \Os^\times := \{g : X \lto \C^\times\}$$
which acts on $\Os$ by multiplication:
$$\G \times \Em \lto \Em, \ \ \ (g, f) \mapsto g f.$$
Since $\dim_\R(X) = 1$, all 2-forms on $X$ are zero.  Hence $\FF_d= \CC_d$.  Identifying $\FF_d$ with $\Omega^1$, there is an  induced action
$$
\G \times \FF_d \lto \FF_d, \ \ \ (g, \nabla) \mapsto \nabla + d\log(g),
$$
where $d\log(g)$ means $g^{-1}dg$.
The set of solutions for $D f = d f = 0$ is precisely $\Os_c$.  In general, the solution set for the equation $\nb f = (d + A)f = 0$ is $\{Ce^{-\int_0^x A} : C \in \Os_c\}$, i.e. $\Os_c$ scaled by $e^{-\int_0^x A}$.  Finally, after applying the gauge $g$, the solution set for the equation
$$g^{-1} \nabla g (f) = (d + A + d \log(g))f = 0$$
is $\{C g(x)^{-1} e^{-\int_0^x A} : C \in \Os_c\}$.
Hence the $\G$-action simply scales the solutions.

On the unit ball $X = (-1, 1)$, every closed 1-form is exact, i.e. $A = d\log(g)$ for some $g \in \G$ or $\nabla$ is equivalent to $d$.  Hence $\Em^\nabla \cong \Os_c$ as an $\Os_c$-module for all connections $\nabla \in \CC_d$.

\subsection{Connections on the unit disk}
Next, let $X$ be the unit disk in $\R^2$ and suppose $f \in \Os$.
Then $f \in \Os_c$ (i.e. constant) if and only if
$$
\frac{\partial f}{\partial x_1} = \frac{\partial f}{\partial x_2} = 0.
$$

Once again, let $r = \rk(\Em) = 1$ for simplicity.  Then $G = \C^\times$ and $\Em \cong \Os$ as an $\Os$-module and $D = d$.  Let $\nb \in \CC_d$.
One would like to say that $f \in \Em$ is a constant if and only if $\nabla f = 0.$  {\em However there is a problem:}  Suppose
$$\nabla f = d f + A f = 0, A \in \Omega^1.$$
Integrating, the solution is of the form
\begin{equation}\label{eq:path_function}
f(z) = C e^{-\int_L A} = f(0)e^{-\int_L A},
\end{equation}
where $L$ is a path connecting $0$ and $z$ in $X$.  Hence for $f$ to be (even locally) well-defined, the integral must depend only on $z$ and independent of $L$.  In other words, $\oint A = 0$ for any small loop that begins and ends at $0$.    By Stokes' theorem, this condition is satisfied if $d A = 0$.  Since $G = \C^\times$ is abelian, for $\nb = d + A$, $$\nb^2 = dA + A \wedge A = d A + \frac{1}{2} [A, A] = d A.$$  This is the motivation for the integrability condition:  If $\nb$ is integrable, then one may locally solve the equation $\nb f = 0$.
Hence $\FF_d$ identifies with $$Z^1 = \{A \in \Omega^1 : d A = 0\}, $$ the space of closed 1-forms when $G$ is abelian.
Meanwhile, the gauge $\G$-action is
$$
\G \times Z^1 \lto Z^1, \ \ (g, \nabla) \mapsto \nabla + d\log(g).
$$
In other words, the gauge $g$ changes the integrable connections by the $\log$-exact form $d \log(g)$.

\subsection{The holomorphic construction}
The integrability condition $\bn ^ 2 = 0$ is trivially satisfied since any $(0,2)$-form is zero in complex dimension $n=1$.  Hence on the unit disk $X$, $\Fh = \Ch$. Moreover, when $r = 1$, every (almost) holomorphic structure $\bn$ is equivalent to $\dbar$ by a gauge transformation.
Similarly, $\Fsf(\nl) = 0$ is also trivially satisfied since all $(2,0)$-forms are zero on $X$.

\section{Interlude: manifolds and functions}\label{sec: smooth manifolds}
We now briefly describe the concept of manifolds.
Let $B \subset \C^n$ (or $B \subset \R^n$) be the open unit ball.  A set $X$ is a manifold if there exist injective maps $\{\phi_i : B \lto X\}_{i \in \II}$, indexed by $\II \subseteq \Z_+$, such that $X = \cup_{i \in \II} \phi_i(B).$  Let $U_i = \phi_i(B)$ and $\UU = \{U_i : i \in \II\}$. For a subset $\JJ \subseteq \II$, let $U_\JJ = \cap_{i \in \JJ} U_i.$
For $i, j \in \II$, define transition maps
$\phi_{ij}   = \phi_i^{-1} \circ \phi_j.$

The injectivity of each $\phi_i$ provides each $U_i$ with the structures on $B$.  However these structures need to be compatible on the intersections $U_{ij}$ in order for $X$ to inherit these structures globally.  These are conditions imposed on the transition maps $\phi_{ij}.$

\begin{defin}
For all $i, j \in \II$,
\begin{enumerate}
\item $X$ is a {\bf topological manifold} if  $\phi_{ij}$ is a homeomorphism.
\item $X$ is a {\bf smooth manifold} if $\phi_{ij}$ is a diffeomorphism.
\item $X$ is a {\bf complex manifold} if $\phi_{ij}$ is holomorphic diffeomorphism.
\end{enumerate}
\end{defin}


Let $X$ be a smooth manifold.  Then each chart $U_i$ inherits an exterior differential operator $d_i$ from $B$.  Let $f : U_i \lto \C$ such that $f \circ \phi_i$ is smooth.
Since each $\phi_{ij}$ is smooth, $f \circ \phi_j$ is smooth for all $j \in \II$.  Together the local exterior differential operators $\{d_i\}_{i \in \II}$ define a global exterior differential operator $d$.

Similarly, suppose that $X$ is a complex manifold.
Then each chart $U_i$ inherits the $\d_i$ and $\db_i$ operators with $d_i = \d_i + \db_i$ from those on $B$.  Let $f : U_i \lto \C$ such that $f \circ \phi_i$ is holomorphic.
Since each $\phi_{ij}$ is a holomorphic diffeomorphism, $f \circ \phi_j$ is holomorphic for all $j \in \II$.  Together these local operators $\{\d_i\}_{i \in \II}$ and $\{\db_i\}_{i \in \II}$ define global differential operators $\d$ and $\db$ and the exterior differential operator decomposes as $d = \d + \dbar$.

\begin{defin}
Let $X$ be a smooth manifold with charts $\UU = \{U_i\}_{i \in \II}$.  Then there is an exterior differential operator $d$ on $X$.  Let $U \subseteq X$ be an open set.  A function $f : U \lto \C$ is {\bf smooth} if $f \circ \phi_i$ is smooth for all $i$.
Suppose that, in addition, $X$ is a complex manifold.  A function $f : U \lto \C$ is {\bf holomorphic} if $f \circ \phi_i$ is holomorphic (i.e $\db_i (f \circ \phi_i) = 0$) for all $i$.
\end{defin}

\begin{rem}\label{rem:abstract manifold}
A smooth manifold $X$ is a topological space locally homeomorphic to $B$ with a globally defined exterior differential operator $d$.  In addition, if $X$ is a complex manifold, then there are globally defined operators $\d$ and $\db$ such that $d = \d + \db$.
\end{rem}
From now on, we shall assume that $X$ is always smooth and that $X$ has a complex structure when discussing holomorphic objects.

\section{The \v{C}ech construction}\label{sec: sheaves}
The most direct way to define (locally) holomorphic and constant functions on $X$ is the \v{C}ech construction.
The manifold $X$ is covered with charts.  We are already familiar with functions on each chart from Sections~\ref{sec:local R-functions} and \ref{sec:local C-functions} where it took only one sentence to define constant and holomorphic functions, respectively, on the unit ball.  The arduous task now is to glue these local functions between charts and into systems of functions.  For this, we introduce sheaves.

Let $X$ be a smooth manifold with charts $\UU = \{U_i\}_{i \in \II}$.
A sheaf on a manifold $X$ assigns to each open set $U \subseteq X$ a particular family of functions.

\subsection{Structure sheaves} For each open set $U \subseteq X$, let
\begin{equation*}
\Os(U)  =  \{f : U \lto \C : f \text{ is smooth }\}.\\
\end{equation*}
In other words, the {\em structure sheaf} $\Os$ on $X$ assigns to each open set $U \subseteq X$ the ring of smooth functions on $U$ (Recall that smoothness makes global sense since each $\phi_{ij}$ is smooth). Let
$$\cC^p(\UU, \Os) = \prod_{i_0 < \cdots < i_p} \Os(U_{i_0 \cdots i_p}), \ \ i_0,  \cdots,  i_p \in \II.$$
The \v{C}ech complex for $(\UU, \Os)$ is
\begin{equation}\label{eq:cech complex}
\cC^\bullet(\UU, \Os) : \ \cC^0(\UU, \Os) \stackrel{d}{\to} \cC^1(\UU, \Os) \stackrel{d}{\to} \cdots  \cC^m(\UU, \Os) \cdots,
\end{equation}
where
$$(d\alpha)_{i_0 \cdots i_{p+1}} = \sum_{j=0}^{p+1} (-1)^j \alpha_{i_0 \cdots \hat{i}_j \cdots i_{p+1}}$$


\subsection{Sheaves of modules}
Local functions on charts $U_i$ and $U_j$ may have different descriptions on $U_{ij}$.  This is described by \v{C}ech 1-cocycles:
\begin{defin}\label{def:cech one cycle}
The {\bf first \v{C}ech $G$-cochain} with respect to $\II$ is
$$
\cC^1(G) = \{\{g_{ij} : U_{ij} \lto G\}_{i,j \in \II}\}.
$$
If the family $\{g_{ij}\}_{i,j \in \II}$ satisfies the additional condition $g_{ij} \circ g_{jk} \circ g_{ki} = \id_{U_{ijk}}$ for all $i,j,k \in \II$, then the family is called a 1-cocycle.  The set of all such 1-cocycles is denoted by $\cF$ (or $\cZ^1(G)$).
\end{defin}
A 1-cocycle $g = \{g_{ij}\}_{i,j \in \II}$ defines a sheaf of $\Os$-module $\Em$ as follows:
Let  $U \subseteq X$ be open.
Define
\begin{equation}\label{eq:sheaf}
\Em(U)  =  \{\{e_i : U \cap U_i \to \C^{\oplus r}\}_{i \in \II} : e_i \text{ is smooth and } g_{ij}.e_j = e_i  \ \  \forall i, j \in \II\}.
\end{equation}
In other words, each element in $\Em(U)$ is a collection of functions such that the function values on $U \cap U_i$ and $U \cap U_j$ differ by the transformation $g_{ij}$.  In this way, $\Em(U)$ is an $\Os(U)$-module.
\begin{rem}\label{rem:sheaf}
Sheaves over $X$ take an open set $U \subseteq X$ and return rings and modules.  Hence the usual algebraic operations such as sum, product, tensor, dualization and etc on sheaves all make sense.  They just mean these operations on the rings and modules on each $U \subseteq X$.

Suppose $\Em$ is a sheaf on $X$.  When referring to an arbitrary open set $U \subseteq X$, we often shorten $\Em(U)$ to $\Em$ when no confusion arises.
\end{rem}
For a sheaf $\Em$, let
$$\cC^p(\UU, \Em) = \prod_{i_0 < \cdots < i_p} \Em(U_{i_0 \cdots i_p}), \ \ i_0,  \cdots,  i_p \in \II.$$
The \v{C}ech complex for $(\UU, \Em)$ is
\begin{equation}\label{eq:cech complexE}
\cC^\bullet(\UU, \Em) : \ \cC^0(\UU, \Em) \stackrel{\nb}{\to} \cC^1(\UU, \Em) \stackrel{\nb}{\to} \cdots  \cC^m(\UU, \Em) \cdots,
\end{equation}
where
$$(\nb\alpha)_{i_0 \cdots i_{p+1}} = \sum_{j=0}^{p+1} (-1)^j \alpha_{i_0 \cdots \hat{i}_j \cdots i_{p+1}}.$$
The resulting cohomologies are denoted by $$\cH^i(\UU, \Em) := \ker(\nb)/\nb(\cC^{i-1}(\UU,\Em)).$$

When $r = 1$ and $g_{ij} = 1 \in \GL(1,\C) = \C^\times$, we recover the structure sheaf $\Em = \Os$.
When $g_{ij} = 1 \in \GL(r,\C)$ is the identity map, the sheaf $\Em \cong \Os^{\oplus r}$ is called trivial.

\subsection{Sheaf isomorphisms} Suppose $\Em^1$ and $\Em^2$ are two sheaves of $\Os$-modules over $X$.  A homomorphism of sheaves $g : \Em^1 \lto \Em^2$ means that for each open $U \subseteq X$, $g_U : \Em^1(U) \lto \Em^2(U)$ is a homomorphism of $\Os$-modules.  In addition if $V \hookrightarrow U$ is open, then
functions on $U$ restrict to functions on $V$, i.e. there is a canonical map $\Em(U) \lto \Em(V)$ for any $\Em$.  We require the following diagram to
commute
$$
\begin{array}{lll}
\Em^1(U) & \stackrel{g_U}{\lto} & \Em^2(U)\\
\downarrow & & \downarrow\\
\Em^1(V) & \stackrel{g_V}{\lto} & \Em^2(V).
\end{array}
$$
The homomorphism $g$ is an isomorphism if $g_U$ is an isomorphism for each open $U \subseteq X$.

Denote by $\cG$ the $G$-zero-cochain
$$
\cG := \cC^0(G) = \{\{g_i : U_i \lto G\}_{i \in \II}\}.
$$
An element in $\cG$ defines a set of isomorphisms $g_i : \Em(U_i) \lto \Em(U_i)$.
\begin{defin}\label{def:Cech groupoid smooth}
The {\bf smooth  \v{C}ech groupoid} $\cF$ consists of
\begin{enumerate}
\item Objects: Sheaves of $\Os$-modules.
\item Equivalence group $\cG$-action:
$$
\cG \times \cF  \lto \cF, \ \ (\{g_i\}_{i \in \II}, \{h_{ij}\}_{i,j \in \II}) \mapsto \{g_i g_j^{-1}h_{ij} \}_{i,j \in \II}.
$$
\end{enumerate}
\end{defin}

\begin{defin}\label{def:Cech groupoid}
 By replacing the word ``smooth'' with ``constant'' starting at this Section~\ref{sec: sheaves}, we obtain the constant structure sheaf $\Os_c$ and the definition of a constant sheaf $\Em_c$ of $\Os_c$-modules and the \v{C}ech complexes $\cC^\bullet(\UU,\Os_c)$ and $\cC^\bullet(\UU,\Em_c)$ and the resulting cohomologies $\cH^i(\UU,\Em_c)$.  If we assume that $X$ is complex in addition, by replacing ``smooth'' with ``holomorphic'', we obtain the holomorphic structure sheaf $\Os_h$ and the definition of a holomorphic sheaf $\Em_h$ of $\Os_h$-modules and the \v{C}ech complexes $\cC^\bullet(\UU,\Os_h)$ and $\cC^\bullet(\UU,\Em_h)$ and the resulting cohomologies $\cH^i(\UU,\Em_h)$.

Denote the resulting locally constant and locally holomorphic \v{C}ech groupoids by $\cF_c$ and $\cFh$, respectively.
\end{defin}
A constant transition function is holomorphic and a holomorphic transition function is smooth.  Hence
\begin{rem}\label{rem:smoothly trivial}
There are natural inclusions
$$\Os_c \hookrightarrow \Os_h \hookrightarrow \Os.$$
A 1-cocycle in $\cFh$ defines a holomorphic sheaf $\Em_h$ but also a smooth sheaf $\Em$ and there is a natural inclusion $\Em_h \hookrightarrow \Em$.  Similarly a 1-cocycle in $\cF_c$ defines a locally constant sheaf $\Em_c$ and also a locally holomorphic sheaf $\Em_h$ and a smooth sheaf $\Em$ resulting in  inclusions $\Em_c \hookrightarrow  \Em_h \hookrightarrow \Em.$  To summarize, there are inclusions
\begin{equation}\label{eq:smoothly trivial}
\begin{array}{ccccc}
\cF_c & \hookrightarrow & \cFh & \hookrightarrow & \cF.
\end{array}
\end{equation}
In the reverse direction,
$$
\Em_h = \Os_h \otimes \Em_c \in \cFh, \ \ \Em = \Os \otimes \Em_h \in \cF, \ \ \Em = \Os \otimes \Em_c \in \cF.
$$
\end{rem}
\noindent {\bf Important assumption:} For the rest of these notes, we shall always assume that the underlying smooth sheaf, that is the final image of the maps (\ref{eq:smoothly trivial}), is isomorphic to the trivial sheaf: $\Os^{\oplus r} \in \cF$.

\section{The smooth de Rham groupoid}\label{sec: global smooth}
The purpose of this section and the next is to describe the sheaves arise from the \v{C}ech constructions in terms of integrable connections, gaining other perhaps more accessible perspectives.
These constructions run in parallel to those of Sections~\ref{sec:local R-functions} and \ref{sec:local C-functions} with much repetition, but this is precisely the point:  The conceptual leap from local construction on unit balls to global construction on manifolds is rather small.  Again the running leitmotif is that a function is locally constant if and only if all its partial derivatives are zero.

\subsection{Structure sheaf}

Let $X$ be a smooth manifold with the exterior differential operator $d : \Os \lto \Omega^1$ (See Remark~\ref{rem:abstract manifold}).  It is then immediate that $\Os_c = \Os^d$ where $\Os^d$ is the sheaf such that for $U \subseteq X$,
$$
\Os^d(U) = \{f \in \Os(U) : df = 0\}.
$$
Again this follows from the principle that $f$ is a constant function if and only if $df = 0$.
\subsection{Tangent, cotangent sheaves and the de Rham complex}
The manifold $X$ is covered by charts in $\UU$ and vector fields are locally defined on each chart and together they form a sheaf over $X$ with transition maps $g_{ij} = {\mathbf D} \phi_{ij} \circ \phi_j^{-1}$, where ${\mathbf D}$ is the standard differential operator on the unit ball.  The resulting smooth sheaf $\T$ is called the tangent sheaf of $X$.
Let $\Om^1$ be the dual sheaf of $\T$ (See Remark~\ref{rem:sheaf}).  It is the sheaf of covector fields on $X$.
By taking exterior products of $\Omega^1$ and extending the exterior differential operators, one obtains the de Rham complex
\begin{equation}\label{eq:smooth de Rham complex}
(\Omega^\bullet, d) :  \Omega^0 \stackrel{d}{\lto} \Omega^{1} \cdots \Omega^i \stackrel{d}{\lto} \Omega^{i+1} \cdots,
\end{equation}
but it is now a complex of sheaves.  If we take the global sections, then
\begin{equation*}
\Omega^0(X) \stackrel{d}{\lto} \Omega^{1}(X) \cdots \Omega^i(X) \stackrel{d}{\lto} \Omega^{i+1}(X) \cdots
\end{equation*}
is an analogue of $\cC^\bullet(\UU,\Os_c)$.

\subsection{Connections}
By Remark~\ref{rem:smoothly trivial}, each constant sheaf $\Em_c$ can be interpreted as a subsheaf of $\Em$.  We begin by fixing a rank-$r$ smooth sheaf of $\Os$-modules:
$$\Em = \{e : X \lto \C^{\oplus r} : e \text{ is smooth }\}.$$
\begin{rem}
The locally constant (and later, holomorphic) sheaves we are considering embed into the trivial sheaf $\Em$ by the assumption at the end of Chapter \ref{sec: sheaves}.
\end{rem}
\begin{defin}\label{def: smooth connection global}
Let $\Omega^i(\Em) := \Omega^i \otimes \Em$.
An $\Os_c$-linear operator
$$\nabla : \Em \lto \Omega^1(\Em)$$
is a {\bf connection} if
$\nabla(f e) = df \otimes e + f \nabla(e)$ for $f \in \Os$ and $e \in \Em$.
Denote by $\CC_d$ the set of all connections on $\Em$.

For each $i$, a connection $\nabla$ extends to an $\Os_c$-linear map
$$
\nabla : \Omega^i(\Em) \lto \Omega^{i+1}(\Em)
$$
satisfying the generalized Leibniz rule
$$\nabla(\eta \otimes e) = d \eta \otimes e + (-1)^i \eta \wedge \nabla(e).$$
\end{defin}
Connections are local operators, operating on $\Em(U)$ for open sets $U \subseteq X$.
As before, the exterior differentiation is generalized to the rank-$r$ trivial connection
$$D : \Em \lto \Omega^1(\Em) := \Omega^1 \otimes \Em, \ \ D(f) = D(\oplus_{i=1}^r f_i) = \oplus_{i=1}^r d f_i.$$
Two connections on $\Em$ differ by an element $A \in \Omega^1(\End(\Em))(X),$ hence,  $\CC_d$ identifies with $\Omega^1(\End(\Em))(X)$.
Then 
$\nabla^2$
is $\Os$-linear, hence, $\nabla^2 \in \Omega^2(\End(\Em))(X)$.
\begin{defin}
The {\bf curvature map} is
$$
\Fsf : \CC_d \lto \Omega^2(\End(\Em)), \ \ \Fsf(\nb) = \nb^2.
$$
For $\nb \in \CC_d$, $\Fsf(\nb)$ is the {\bf curvature} of $\nb$.
A connection $\nabla$ is {\bf integrable} if $\Fsf(\nabla) = 0$.
If $\nabla$ is integrable, then $$\Em^\nabla(U) = \{e \in \Em(U) : \nabla e = 0\} , \ \ \forall U \subseteq X$$ is the locally constant sheaf of $\Os^d$-modules with respect to $\nabla$.
Denote by $\FF_d$ (or $\ZZ^1(G)(X)$) the subspace of integrable connections in $\CC_d$.
\end{defin}


When $\nb$ is integrable, the de Rham complex is
$$
(\Omega^{\bullet}(\Em), \nb) :  \Omega^{0}(\Em) \stackrel{\nb}{\lto} \Omega^{1}(\Em) \cdots \Omega^{i}(\Em) \stackrel{\nb}{\lto} \Omega^{i+1}(\Em) \cdots.
$$
Taking global sections
$$
\Omega^{0}(\Em)(X) \stackrel{\nb}{\lto} \Omega^{1}(\Em)(X) \cdots \Omega^{i}(\Em)(X) \stackrel{\nb}{\lto} \Omega^{i+1}(\Em)(X) \cdots,
$$
the resulting $i$-th cohomology is $$\H_\nb^i(\Em) = \ker(\nb)/\nb(\Omega^{i-1}(\Em)(X)).$$  These are the analogues of $\cC^\bullet(\UU,\Em_c)$ and $\cH^i(\UU,\Em_c)$.

The sheaf $\T$ generates the sheaf $\D$ of rings of differential operators on $\Os$.  Integrability of $\nb$ is precisely the condition for $\Em$ to inherit a $\D$-module structure from its $\T$-module structure.

\subsection{The gauge group and its action}
The group $G$ acts on $\C^{\oplus r}$ by matrix-vector multiplication.  Define the $G$-gauge group
$$\G = \Omega^0(G)(X) := \Os(G)(X) = \{g : X \lto G\}.$$
Then $\G$ acts on $\Em$ by $\Os$-module automorphism
$$
\G \times \Em \lto \Em, \ \ (g,e) \mapsto ge, \ \ ge(x) = g(x)e(x),
$$
where $g(x)e(x)$ means vector-matrix multiplication at $x \in X$.
Suppose $g \in \G$.  Then $g : \Em \stackrel{\cong}{\lto} \Em$.  If $\nb$ is a connection on $\Em$, then $g$ pulls back $\nb$ to a connection $g^*\nb$ on $\Em$ as follows:
$$
g^*\nb : \Em \stackrel{g}{\lto} \Em \stackrel{\nb}{\lto} \Omega^1(\Em) \stackrel{g^{-1}}{\lto} \Omega^1(\Em).
$$
This give a $\G$-action on $\CC_d$.
Identifying $\CC_d$ with $\Omega^1(\End(\Em))(X)$, $\nb = D + A$ and the $\G$-action is
$$
\G \times \CC_d \lto \CC_d, \ \ \ (g,\nabla) \mapsto g^{-1} \nabla g = D + g^{-1} (Dg) + g^{-1} A g.
$$

The curvature of $g^{-1} \nb g$ is $(g^{-1} \nb g)^2 = g^{-1} \nb^2 g.$
Hence $\nb^2 = 0$ if and only if $g^{-1}\nb^2g = 0$ and the gauge action restricts to an action
$\G \times \FF_d \lto \FF_d.$
\begin{defin}\label{def: smooth de Rham groupoid}
The {\bf de Rham groupoid} $\FF_d$ is
\begin{enumerate}
\item Objects: the family of integrable connections on $\Em$.
\item Equivalence: gauge group $\G$-action
$$
\G \times \FF_d \lto \FF_d, \ \ (g, \nabla) \mapsto g^{-1} \nabla g.
$$
\end{enumerate}
\end{defin}
The connection $D$ defines a map
$$
\Dsf : \G \lto \CC_d, \ \ \Dsf(g) = g^{-1} D g.
$$
There is a sequence of maps
$$
\Omega^0(G)(X) \stackrel{\Dsf}{\lto} \Omega^1(\End(\Em))(X) \stackrel{\Fsf}{\lto} \Omega^2(\End(\Em))(X).
$$
We have already seen in Section~\ref{sec:local R-functions} that over a small enough open ball $U \subseteq X$, this sequence is exact.  However on large open sets and especially over $X$, this may not be the case and the resulting orbit space is our main focus.

\begin{rem}\label{rem:global smooth E_c = E^nb}
For an integrable connection $\nb \in \FF_d$ and a small enough $U \subseteq X$, $\Em^\nb(U)$ is isomorphic to $\Os_c^{\oplus r}(U)$ as an $\Os_c(U)$-module as stated in Remark~\ref{rem: local consant subset smooth}.  Hence $\Em^\nb$ is isomorphic to a locally constant sheaf $\Em_c \in \cF_c$ on $X$.
\end{rem}

\section{The holomorphic de Rham and Dolbeault groupoid}\label{sec:global holomorphic}

Let $X$ be a complex manifold with operators $d = \d + \dbar$ (See Remark~\ref{rem:abstract manifold}).  There is a decomposition of the sheaf of exterior $m$-forms according to types and projections:
\begin{equation}\label{eq:hodge decomposition}
\Omega^m = \bigoplus_{q+p=m}\Omega^{q,p}, \ \ \ P_{q,p} : \Omega^m \lto \Omega^{q,p}.
\end{equation}
For $d : \Omega^m \lto \Omega^{m+1}$,
$$
d = \d + \db, \ \ \ \d = P_{q+1,p} \circ d \ \ \text{ and } \ \ \ \db = P_{q,p+1} \circ d.
$$
Moreover, for $p$ and $q$, the de Rham complex decomposes as
\begin{eqnarray}\label{eq: Dolbeault complex}
(\Omega^{q,\bullet}, \dbar) :  \Omega^{q,0} \stackrel{\dbar}{\lto} \Omega^{q,1} \cdots \Omega^{q,i} \stackrel{\dbar}{\lto} \Omega^{q,i+1} \cdots,\\
(\Omega^{\bullet,p}, \d) :  \Omega^{0,p} \stackrel{\d}{\lto} \Omega^{1,p} \cdots \Omega^{i,p} \stackrel{\d}{\lto} \Omega^{i+1,p} \cdots.
\end{eqnarray}

The complex structure on $X$ may be equivalently described by an $\Os$-linear operator $J : \T \lto \T$ with $J^2 = -1$ such that the Nijenhuis tensor of $J$ vanishes, i.e.
\begin{equation}\label{eq:NJ}
N_J(X,Y) := [X,Y] + J([JX,Y] + [X,JY]) - [JX, JY] = 0
\end{equation}
for all $X, Y \in \T$.
This is the content of the Newlander-Nirenberg theorem \cite{NN1}.  The smooth tangent sheaf decomposes as
$\T = \T^{1,0} \oplus \T^{0,1}$ by the eigenspaces of $J$ with eigenvalues ${\bf i}$ and $-{\bf i}$, respectively.  This induces a decomposition $\Omega^1 = \Omega^{1,0} \oplus \Omega^{0,1}$ which gives rise to the decomposition (\ref{eq:hodge decomposition}) for each $m$.  This again gives a decomposition $d = \d + \db$ which is equivalent to $\db^2 = 0$ and to Equation (\ref{eq:NJ}).

The sheaf $\Os_h$ is equivalent to $\Os^\dbar$, where
\begin{equation}\label{eq:holomorphic structure sheaf}
\Os^\dbar(U) = \{f \in \Os(U) : \dbar f = 0\}.
\end{equation}
The holomorphic tangent sheaf is
$$\T_h = \{V \in \T : V(f) \in \Os_h,  \forall f \in \Os_h\}.$$
The holomorphic cotangent sheaf $\Omega_h^1$ is the $\Os_h$-dual of $\T_h$ and $\Omh^q = \wedge_{i=1}^q \Omh^1$.  The holomorphic de Rham complex is
\begin{equation}\label{eq:holomorphic de Rham complex}
(\Omh^{\bullet}, \d) :  \Omh^{0} \stackrel{\d}{\lto} \Omh^{1} \cdots \Omh^{i} \stackrel{\d}{\lto} \Omh^{i+1} \cdots.
\end{equation}

\subsection{Holomorphic structures}
\begin{defin}
Let $\Omega^{q,p}(\Em) := \Omega^{q,p} \otimes \Em$.
An almost holomorphic structure on $\Em$ is an operator
$$\bn : \Em \lto \Omega^{0,1}(\Em), \ \ \text{ with } \bn(f e) = \dbar f \otimes e + f \bn(e)$$
for $f \in \Os$ and $e \in \Em$.  Denote by $\Ch$ the space of almost holomorphic structures on $\Em$.

For $p+q=i$, $\bn$ extends to an $\Os_c$-linear map
$$
\bn : \Omega^{p,q}(\Em) \lto \Omega^{p,q+1}(\Em)
$$
satisfying the generalized Leibniz rule
$$\bn(\eta \otimes e) = \dbar \eta \otimes e + (-1)^i \eta \wedge \bn(e).$$
\end{defin}

The operator $D$ also decomposes according to types $ D = D' + D'':$
$$D' : \Em \lto \Omega^{1,0}(\Em) := \Omega^{1,0} \otimes \Em, \ \ D'(f) = D'(\oplus_{i=1}^r f_i) = \oplus_{i=1}^r \d f_i.$$
$$D'' : \Em \lto \Omega^{0,1}(\Em) := \Omega^{0,1} \otimes \Em, \ \ D''(f) = D(\oplus_{i=1}^r f_i) = \oplus_{i=1}^r \dbar f_i.$$
Define map
$$
\Fsf : \Ch \lto \Omega^{0,2}(\End(\Em)), \ \ \Fsf(\bn) = \bn^2.
$$
\begin{defin}\label{def:holomorphic structure}
A {\bf holomorphic structure} $\bn$ on $\Em$ is an almost holomorphic structure satisfying the integrability condition $\Fsf(\bn) = 0$.
Denote by $\Fh$ (or $\ZZ^{0,1}(G)(X)$) the holomorphic structures on $\Em$.
\end{defin}
Then $D'' \in \Fh$ and
$$
\Em^{D''} = \{e \in \Em : D'' e = 0\} \cong \Os_h^{\oplus r}
$$
is a sheaf of $\Os^\dbar$-module.

The difference of two almost holomorphic structure is $\Os$-linear; hence, the space of almost holomorphic structures $\Ch$ identifies with $\Omega^{0,1}(\End(\Em))(X)$.
If $\bn \in \Fh$,
then $\Em^\bn = \{f : \bn f = 0\}$ is the sheaf of $\Os_h$-modules of $\bn$-holomorphic sections and the $\bn$-Dolbeault complex is
$$
(\Omega^{q,\bullet}(\Em), \bn) :  \Omega^{q,0}(\Em) \stackrel{\bn}{\lto} \Omega^{q,1}(\Em) \cdots \Omega^{q,i}(\Em) \stackrel{\bn}{\lto} \Omega^{q,i+1}(\Em) \cdots.
$$
Take global sections
$$
\Omega^{q,0}(\Em)(X) \stackrel{\bn}{\lto} \Omega^{q,1}(\Em)(X) \cdots \Omega^{q,i}(\Em)(X) \stackrel{\bn}{\lto} \Omega^{q,i+1}(\Em)(X) \cdots.
$$
The resulting $i$-th cohomology is $$\H_\bn^{q,i}(\Em) = \ker(\bn)/\bn(\Omega^{q,i-1}(\Em)(X)).$$  These are the analogues of $\cC^\bullet(\UU,\Omega_h^q(\Em_h))$ and $\cH^i(\UU,\Omega_h^q(\Em_h))$.

\subsection{The gauge group action I}
The gauge group $\G$ acts on the almost holomorphic structures
$$
\G \times \Ch  \lto \Ch, \ \ (g, \bn) \mapsto g^{-1} \bn g.
$$
This $\G$-action preserves $\Fh$.
\begin{defin}\label{def:holomorphic groupiod}
The {\bf groupoid of holomorphic structures} $\Fh$ is
\begin{enumerate}
\item Objects: set of holomorphic structures on $\Em$.
\item Equivalence: gauge group $\G$-action:
$$
\G \times \Fh \lto \Fh, \ \ (g, \bn) \mapsto g^{-1} \bn g.
$$
\end{enumerate}
\end{defin}
The operator $D''$ defines a map
$$
\Dsf : \G \lto \Ch, \ \ \Dsf(g) = g^{-1}D''g
$$
Identifying $\Ch$ with $\Omega^{0,1}(\End(\Em))(X)$, there is a sequence of maps
$$
\G \stackrel{\Dsf}{\lto}  \Omega^{0,1}(\End(\Em))(X) \stackrel{\Fsf}{\lto} \Omega^{0,2}(\End(\Em))(X).
$$
\begin{rem}\label{rem:global holomorphic E_h = E^bn}
For a holomorphic structure $\bn \in \Fh$ and a small enough $U \subseteq X$, $\Em^\bn(U)$ is isomorphic to $\Os_h^{\oplus r}(U)$ as an $\Os_h(U)$-module as stated in Remark~\ref{rem: local holomorphic subset smooth}.  Hence $\Em^\bn$ is isomorphic to a holomorphic sheaf $\Em_h \in \cFh$ on $X$.
\end{rem}

\subsection{Holomorphic $\lambda$-connections}
A local function $f \in \Os(U)$ is constant if and only if $\d f = \dbar f = 0$.  Hence,
\begin{equation}\label{eq:constant structure sheaf}
\Os_c(U) = \Os_h^\d(U) = (\Os^\dbar)^\d(U) = \{f \in \Os(U) : \dbar f = \d f = 0\}.
\end{equation}
Fix a holomorphic structure $\bn$.  By Remark~\ref{rem:global holomorphic E_h = E^bn}, $\Em^\bn \cong \Em_h \in \cFh$ and we will use $\Em_h$ to simplify notations.  We emphasize here that $\Em_h$ may no longer be isomorphic to $\Os_h^{\oplus r}$.  Let $\lambda \in \C$.
\begin{defin}
Let
$\Omh^{i}(\Em_h) := \Omh^i \otimes \Em_h$ and
$\lambda \in \C$.  A {\bf holomorphic $\lambda$-connection} or {\bf $\lambda$-connection} on $\Em_h$ is an operator $$\nl : \Em_h \lto \Omh^{1}(\Em_h)$$
satisfying $\nl(f e) = \lambda \d f \otimes e + f \nl(e)$ for $f \in \Os_h$ and $e \in \Em_h$.  Denote by $\Chl$ the space of all $\lambda$-connections on $\Em_h$.

For each $i$, $\nl$ extends to an $\Os_c$-linear map
$$
\nl : \Omh^i(\Em) \lto \Omh^{i+1}(\Em)
$$
satisfying the generalized Leibniz rule
$$\nl(\eta \otimes e) = \lambda \d \eta \otimes e + (-1)^i \eta \wedge \nl(e).$$
\end{defin}
The curvature map is
$$
\Fsf : \Chl \lto \Omega_h^2(\End(\Em)), \ \ \Fsf(\nl) = (\nl)^2.
$$
\begin{defin}\label{def:l-connection}
A $\lambda$-connection $\nl \in \Chl$ is {\bf integrable} if $\Fsf(\nl) = 0$.
Denote by $\Fhl$ the space of integrable $\lambda$-connections on $\Em_h$.
\end{defin}
Notice that $\lambda D'$ is an integrable $\lambda$-connection on $\Os_h^{\oplus r}$.  In general, there is no such canonical $\lambda$-connection on $\Em_h$.
If $\nl, \nl_0 \in \Chl$, then
$\nl - \nl_0 \in \Omega_h^1(\End(\Em_h))(X).$
Hence the space $\Chl$ of integrable $\lambda$-connections identifies with $\Omega_h^1(\End(\Em_h))(X)$.
If $\nl$ is integrable, then $\Em_h^{\nl} = \{e : \nl e = 0\}$ is the sheaf of $\Os_c$-modules of constant sections with respect to $\nl$ and the $\nl$-de Rham complex is
$$
(\Omh^{\bullet}(\Em_h), \nl) :  \Omh^{0}(\Em_h) \stackrel{\nl}{\lto} \Omh^{1}(\Em_h) \cdots \Omh^{i}(\Em_h) \stackrel{\nl}{\lto} \Omh^{i+1}(\Em_h) \cdots.
$$

\begin{rem}
Vector fields in $\T_h$ are first order differential operators on $\Os_h$.  By composing these operators, $\T_h$ generates $\D_h$, the sheaf of rings of holomorphic differential operators.  Integrability of $\nb^1$ is precisely the needed condition for $\Em_h$ to inherit a $\D_h$-module structure from its $\T_h$-module structure \cite{Bo1, Co2, Sc1}.
\end{rem}



\begin{defin}\label{def:holomorphic de Rham groupoid}
The groupoid $\Fl$ on $\Em$ consists of
\begin{enumerate}
\item Objects: $\Fl = \{(\bn,\nbl) : \bn \in \Fh, \nbl \in \Fhl, \text{ where } \Em_h = \Em^\bn\}$.
\item Equivalence: gauge group $\G$:
$$
\G \times \Fl \lto \Fl, \ \ (g, (\bn, \nbl)) = (g^{-1}\bn g, g^{-1} \nbl g).
$$
\end{enumerate}
$\Fl$ is called the {\bf holomorphic de Rham groupoid} when $\lambda = 1$ and the {\bf Dolbeault groupoid} when $\lambda = 0$.
\end{defin}
%

\begin{rem}\label{rem:F_d = F_h_d}
Let $\lambda \in \C^\times$.
For a small enough open set $U \subseteq X$,
$$\Em^\nb(U) \cong (\Os^{\oplus r})^D(U) \cong (\Os_h^{\oplus r})^{\lambda D''}(U) \cong (\Em^\bn)^{\nbl}(U)$$
by Remarks~\ref{rem: local holomorphic subset smooth} and \ref{rem: local consant subset holomorphic subset smooth}.  Hence, $\FF_d$ and $\Fl$ are equivalent constructions.
\end{rem}

\begin{rem}\label{rem:underlying holomorphic structure}
There is a natural map (functor)
\begin{equation*}
\Fl \lto \Fh, \ \ (\bn, \nb) \mapsto \bn.
\end{equation*}
\end{rem}
\begin{rem}\label{rem: Higgs embedding}
When $\lambda = 0$, $\nb^0$ is $\Os_h$-linear and $(\Em_h, \nb^0)$ is called a Higgs bundle.  Moreover there is a natural embedding
$$\Fh \lto \FF_d^0, \ \ \Em_h \mapsto (\Em_h, 0).$$
This is not true when $\lambda \neq 0$ since there is no canonical integrable $\lambda$-connection on a general $\Em_h$.
\end{rem}
The holomorphic gauge group with respect to $\Em_h$ is
$$
\G_h = \{g \in \G : g : \Em_h \to \Em_h\}.
$$
There is an action
$\Gh \times \Em_h \lto \Em_h$.
For fixed $\bn \in \Fh$ and $\lambda \in \C^\times$, define
\begin{defin}\label{def:holomorphic de Rham groupoid fixing h}
The groupoid $\Fhl(\Emh)$ consists of
\begin{enumerate}
\item Objects: the set of integrable $\lambda$-connections on $\Emh$.
\item Equivalence: $\Emh$-holomorphic gauge group $\Gh$-action:
$$
\Gh \times \Fhl(\Emh)  \lto \Fhl(\Emh), \ \ (g, \nbl) \mapsto g^{-1} \nbl g.
$$
\end{enumerate}
\end{defin}


\begin{rem}\label{rem: global lambda connection}
Fix $\Em_h \in \Fh$ and let $\lambda \in \C^\times$.  If $\nb^{1} \in \FF_h^1(\Emh)$, then $\lambda \nb^1 \in \Fhl(\Emh)$.
This is an equivalence of $\Fh^1(\Emh)$ and $\Fhl(\Emh)$.
\end{rem}

\section{The equivalence of groupoids}\label{sec:iso}
A functor $\FFF : \CC_1 \lto \CC_2$ between groupoids takes objects to objects and equivalence relations to equivalence relations, satisfying the general axioms concerning categories and functors.  A functor is an equivalence of groupoid if it is fully faithful and surjective onto the isomorphism classes \cite{GM1, GX1}.  In the previous sections, we described how objects in $\FF_d, \Fl$ ($\lambda \in \C^\times$) give rise to objects in $\cF_c$ (Remark~\ref{rem:global smooth E_c = E^nb} and \ref{rem:global holomorphic E_h = E^bn}).  In other words, we have described functors $\FF_d \lto \cF_c$ and $\Fl \lto \cF_c$.  In this section, we bring in the groupoid $\Hom(\pi,G)$.


Let $X$ be a manifold and $\pi = \pi_1(X,x)$ the fundamental group of $X$ as in Section~\ref{sec:Hom} and $p : \tX \lto X$ the universal cover.  Let $\tOs$ be the structure sheaf of $\tX$, $\tEm = \tOs^{\oplus r}$ and $\td, \tD$ be the exterior differentiation and the trivial connection on $\tEm$.    As before, $\tEm^{\tD} \cong \tOs_c^{\oplus r}$.  The projection $p$ induces natural maps $p^* : \Os \lto \tOs, \ \ p^* : \Os_c \lto \tOs_c$.  The images are $\pi$-invariant local functions.


Suppose $\Em_c \in \cF_c$ is a locally constant sheaf.  Let $\gamma$ be a loop representing an element $[\gamma] \in \pi$.  Starting at $x$, there is a sequence $\{i_1, \cdots, i_k\} \in \II$ with $U_{i_1}, \cdots U_{i_k}$ following and covering $\gamma$.  Since the transition maps $g_{i_j, i_{j+1}}$ are constant maps, we obtain an element $g = \prod g_{i_j, i_{j+1}}$ for $[\gamma]$.  The 1-cocycle condition (See Definition~\ref{def:cech one cycle}) guarantees that $g$ depends only on $[\gamma]$.  This gives a representation of the fundamental group $\rho : \pi \lto G$.

Let $\rho \in \Hom(\pi,G)$.
Then $\pi$ acts on $\tX$ freely by deck transformation and on $\C^r$ by the linear $G$-action via $\rho$.   Let $\{e_i : X \lto \C^{\oplus r}\}_{i=1}^r$ be a basis and $\{\te_i : \tX \lto \C^{\oplus r}\}_{i=1}^r$ be an equivariant basis with respect to $\rho$, i.e.
$$\rho(\sigma).\te_i(x) = \te_i(\sigma . x).$$
Let $U \subseteq X$, $\tU = p^{-1}(U) \subseteq \tX$ and $f \in \Em(U)$.  Then $f = \sum_{i=1}^r f_i e_i$ where $f_i : U \lto \C \text{ for } 1\le i \le r.$  Let $$\tf : \tU \lto \C^{\oplus r}, \ \ \tf(\ty) = \sum_{i=1}^r p^*(f_i)(\ty) \te_i(\ty).$$
Then $\tD\tf \in \Omega^1(\tEm)(\tU)$ is equivariant with respect to the $\pi$-action, hence, descends to an element $\eta \in \Omega^1(\Em)(U)$.  This defines a connection $\nabla : \Em \lto \Omega^1(\Em).$   It is integrable since $\tD$ is integrable. Notice that the construction of $\nabla$ involves a choice of basis and is unique only up to gauge transformations.  To be specific, let
$$
\tD \te_i = \sum_{j = 1}^r a_{ij} \te_j, \ \ a_{ij} \in \Omega^1(\tOs), \text{ i.e. } a_{ij} \text{ is a 1-form on } \tX.
$$
Let $\tA$ be the $r \times r$ matrix with entries $a_{ij}$.  Then $\tA \in \Omega^1(\tEm)$ is $\pi$-invariant.  Hence $\tA = p^*(A)$ for some $A \in \Omega^1(\Em)$.  Then $\nb = D + A$ is the desired integrable connection whose monodromy representation is in $[\rho]$.

By Remark~\ref{rem:global smooth E_c = E^nb},
each $\Em^\nb \in \FF_d$ is isomorphic to a locally constant sheaf $\Em_c \in \cF_c$.
By Remark~\ref{rem:F_d = F_h_d},
$\FF_d$, $\FF_d^\lambda$ are equivalent constructions for $\lambda \in \C^\times$.

To summarize, let $\lambda \in \C^\times$.  Then there are the following correspondences of groupoids:
\begin{equation}\label{eq:equivalence of groupoids 1}
\begin{array}{ccccccccccc}
& & \Hom(\pi,G) & &\\
& \nearrow & & \searrow  \\
\cF_c & & \longleftarrow & & \FF_d& & \longleftrightarrow & & \FF_d^\lambda
\end{array}
\end{equation}
In fact, these groupoids are all equivalent.   We refer readers to \cite{GM1, GX1} for a proof.  What we have done in effect is to show that these functors are surjective on isomorphism classes (one can go around the triangle and see that objects that one begins and ends with are equivalent).

\subsection{K\"ahler manifolds and projective varieties}
Let $\P^n$ be the projective $n$-space over $\C$.  Finite sets of homogeneous polynomials in $(n+1)$ variables define algebraic varieties (loci) in $\P^n$.
\begin{defin}
A {\bf smooth projective variety} $X$ of dimension $n$ over $\C$ is a projective variety that is also a manifold.
\end{defin}
Smooth projective varieties over $\C$ are compact K\"ahler manifolds.  However compact K\"ahler manifolds are not necessarily projective varieties.
If $X$ is a compact K\"ahler manifold, then Equivalence~(\ref{eq:equivalence of groupoids 1}) extends to the case of $\lambda = 0$.  We do not provide a proof here and refer readers to \cite{Co1, Do1, GM1, GX1, Hi1, Si0}, but emphasize that this last equivalence is complicated.  For example, for the correspondence $(\Em_h^1,\nb^1) \mapsto (\Em_h^2, \nb^0)$, it is not the case that $\Em_h^1 \cong \Em_h^2$ in general.

\section{Moduli spaces by simple examples}\label{sec:simple examples}
From the groupoids $\Hom(\pi,G), \cF_c, \cFh, \FF_d, \Fh$ and $\Fl$ arise the moduli spaces $\MB, \cM_c, \cM_h, \MdR, \M_h$ and $\MdR^\lambda$, universal geometric objects parameterizing their respective isomorphism classes as outlined in Section~\ref{sec:groupoid}.  The word universal here means that the resulting moduli spaces represent (or corepresent) certain functors.  This is firmly in the realm of Geometric Invariant Theory or GIT for short.  We refer readers to \cite{Mu1} for a full dosage and \cite{Ne1} for a gentler introduction.  The subject has been continuously exploding for years and we point to the reference section of \cite{Mu1} as evidence.
The process of forming moduli spaces is a subtle and delicate one as demonstrated in Section~\ref{sec:Hom} and the resulting moduli space may not be able to encompass the entire family of isomorphism classes.

The Betti moduli space $\MB$ is relatively simple to construct and we have already done that in Section~\ref{sec:Hom}; however, recall the subtle problems when forming quotients.  This space has an algebraic, analytic and smooth structure from that of $G$.


Here we
give a few simple examples of these moduli spaces when $G$ is the abelian $\GL(1,\C) = \C^\times$.  This puts us back in the more familiar situation of the usual first cohomology, where $\Em$, $\Em_h$ and $\Em_c$ are of rank-1.  Then $\End(\Em) = \Os$, $\End(\Em_h) = \Os_h$ and $\End(\Em_c) = \Os_c$.

\subsection{Classical and abelian cohomology}
Our theory is developed for reductive groups $G$ and is an extension of the theory of the Picard variety where $G = \C^\times$.  This in turn arises from an understanding of the classical cohomology theory where $G$ is replaced by the additive group $\C$ which is not reductive.

\subsubsection{The case of the additive group $\C$}  This is the setting of classical cohomologies.
The constant and holomorphic complexes compute the cohomologies $\cH^i(\UU,\Os_c)$ and $\cH^p(\UU, \Omega_h^q)$:
\begin{eqnarray}
\cC^\bullet(\UU, \Os_c) & : &  \cC^0(\UU, \Os_c) \stackrel{d}{\to} \cC^1(\UU, \Os_c) \stackrel{d}{\to} \cdots  \cC^j(\UU, \Os_c) \cdots,\\
\cC^\bullet(\UU, \Omega_h^q) & : &  \cC^0(\UU, \Omega_h^q) \stackrel{d}{\to} \cC^1(\UU, \Omega_h^q) \stackrel{d}{\to} \cdots  \cC^j(\UU, \Omega_h^q) \cdots.
\end{eqnarray}
There is the notion of {\em acyclic cover} (\S 3, Chapter 0, \cite{GH1}).  These acyclic covers exist under mild assumptions.
When $\UU$ is an acyclic cover with respect to $\Omega_h^q$, $\cC^\bullet(\UU, \Omega_h^q)$ computes the \v{C}ech cohomologies $\cH^p(\Omega_h^q)$ by definition.  When $\UU$ is an acyclic cover with respect to $\Os_c$, $\cC^\bullet(\UU, \Os_c)$ computes the \v{C}ech cohomologies $\cH^i(\Os_c)$ by definition.  The fact that \v{C}ech cohomologies are independent of acyclic covers is known as the Leray theorem.

Given a sheaf or a complex of sheaves $F^\bullet$ on $X$, there is a notion of acyclic resolution (\S 1, Chapter 3, \cite{Ha1}).  Let $F^\bullet \lto I^\bullet$ be such an acyclic resolution.  Then the sheaf cohomologies $F^\bullet$ are computed via $I^\bullet$ after taking global sections.   As an example, when $\UU$ is an acyclic cover, $\cC^\bullet(\UU, \Os_c)$ and $\cC^\bullet(\UU, \Omega_h^q)$ should be thought of as the results of applying the global section functor to the sheafified acyclic \v{C}ech complex resolutions (\S 4, Chapter 3, \cite{Ha1}). A fundamental result of homological algebra is that the resulting cohomological groups are independent of the acyclic resolutions.   

The de Rham and Dolbeault complexes are respectively:
\begin{eqnarray}
(\Omega^\bullet, d) & : & \Omega^0 \stackrel{d}{\lto} \Omega^{1} \cdots \Omega^j \stackrel{d}{\lto} \Omega^{j+1} \cdots, \\
(\Omega^{q,\bullet}, \dbar) & : & \Omega^{q,0} \stackrel{\dbar}{\lto} \Omega^{q,1} \cdots \Omega^{q,j} \stackrel{\dbar}{\lto} \Omega^{q,j+1} \cdots
\end{eqnarray}
The Dolbeault complex is an acyclic resolution of the holomorphic sheaf $\Omega_h^q$, hence, computes the holomorphic \v{C}ech cohomology:
$$\H_\dbar^{q,p}(\Os) \cong \cH^p(\Omega_h^q).$$
This is the content of Dolbeault theorem.

The de Rham complex is an acyclic resolution of the constant sheaf $\Os_c = \C$, hence, computes the constant \v{C}ech cohomology:  $$\H_d^i(\Os) \cong \cH^i(\Os_c).$$
This is the content of de Rham theorem.

The de Rham complex also resolves the holomorphic de Rham complex
\begin{equation}\label{eq:holomorphic de Rham}
\begin{array}{ccccccccccc}
(\Omega^\bullet, d) & : & \Omega^0 & \stackrel{d}{\lto} & \Omega^{1} & \cdots & \Omega^j & \stackrel{d}{\lto} & \Omega^{j+1} & \cdots\\
& & \uparrow & & \uparrow & &  \uparrow & & \uparrow \\
(\Omh^{\bullet}, \d) & : & \Omh^{0} & \stackrel{\d}{\lto} & \Omh^{1} & \cdots & \Omh^{j} & \stackrel{\d}{\lto} & \Omh^{j+1} & \cdots.
\end{array}
\end{equation}
The decomposition (\ref{eq:hodge decomposition}) gives a filtration of the de Rham complex
$$F^k\Omega^i = \bigoplus_{q+p = i, q \ge k} \Omega^{q,p}.$$
This then induces a filtration on the de Rham cohomology which gives rise to a double complex and a spectral sequence with $E_1^{0,1} = \H_\dbar^{0,1}(\Os)$.  The projection to the graded piece gives a map
\begin{equation}\label{eq:projection to the graded}
\H_d^1(\Os) \lto E_1^{0,1} = \H_\dbar^{0,1}(\Os).
\end{equation}
If $X$ is smooth projective, then the spectral sequence degenerates at $E_1$ and the Hodge theorem states that the de Rham cohomologies decompose as
$$
\H_d^i(\Os) \cong \bigoplus_{q+p = i} \H_\dbar^{q,p}(\Os) \cong \bigoplus_{q+p = i} \cH^p(\Omega_h^q).
$$

By covering based loops with open charts as in Section~\ref{sec:iso}, there is an isomorphism $\cH^1(\Os_c) \cong \Hom(\pi,\C) $.  To summarize, we have the following equivalences:
\begin{equation}\label{eq:equivalence of groupoid example}
\begin{array}{ccccccccccc}
\Hom(\pi,\C) & \cong & \cH^1(\Os_c) & \cong & \H_d^1(\Os) & \stackrel{X \text{ projective }}{\cong} & \H^1(\Omega_h^0) \oplus \H^0(\Omega_h^1).\\
\end{array}
\end{equation}

\subsubsection{$G = \C^\times$}  The moduli spaces are generalizations of classical cohomologies.
Let $\Os_c^\times$ (resp. $\Os_h^\times$) be the sheaf of (multiplicative) locally constant (resp. holomorphic) functions that are never zero.  These sheaves are related by exact sequences

\begin{eqnarray} \label{eq:exp c sequences}
0 \lto \Z \lto \Os_c \stackrel{exp}{\lto} \Os_c^\times \lto 0,\\ \label{eq:exp h sequences}
0 \lto \Z \lto \Os_h \stackrel{exp}{\lto} \Os_h^\times \lto 0.
\end{eqnarray}
Here $\Z$ denotes the sheaf of constant functions with integer values.
The long cohomological sequence associated with the short exact sequence (\ref{eq:exp c sequences}) is
\begin{equation}\label{eq:constant long exact sequence}
\cdots \stackrel{0}{\lto} \cH^1(\Z) \lto \cH^1(\Os_c) \lto \cH^1(\Os_c^\times) \stackrel{\delta}{\lto} \cH^2(\Z) \cdots.
\end{equation}
Then $\cM_c := \ker(\delta) \cong \cH^1(\Os_c)/\cH^1(\Z)$.  Here $\delta$ is the Chern class. 
Replacing the subscript $c$ with $h$ in the above construction, the long exact sequence associated with sequence (\ref{eq:exp h sequences}) is
\begin{equation}\label{eq holomorphic long exact sequence}
\cdots \cH^1(\Z) \lto \cH^1(\Os_h) \lto \cH^1(\Os_h^\times) \stackrel{\delta}{\lto} \cH^2(\Z) \cdots.
\end{equation}
The cohomology $\cH^1(\Os_h^\times)$ is the moduli space of all rank-1 holomorphic sheaves on $X$ and $\cM_h := \ker(\delta) \cong \cH^1(\Os_h)/\cH^1(\Z)$ which is the moduli space of rank-1 sheaves that embed into the trivial smooth sheaf (See Remark~\ref{rem:smoothly trivial}).

\subsection{The punctured disk}\label{subsec: punctured disk}
Let $X$ be the punctured unit disk in $\C$ and let $z = x + i y$.

\subsubsection{The Betti moduli space} The fundamental group $\pi$ is isomorphic to $\Z$ and $\Hom(\pi_1(X), \C^\times) = \C^\times$.

\subsubsection{The \v{C}ech moduli}
From the long exact sequence (\ref{eq:constant long exact sequence}), the moduli space of rank-1 locally constant sheaves is
$$\cM_c = \ker(\delta) \cong \cH^1(\Os_c)/\cH^1(\Z) \cong \C/\Z$$
which is a 1-dimension complex torus as is $\MB$.

From the long exact sequence (\ref{eq holomorphic long exact sequence}),
the moduli space of all rank-1 holomorphic sheaves on $X$ is $\cM_h = \ker(\delta)$ which is the moduli space of rank-1 sheaves that embed into the trivial smooth sheaf (See Remark~\ref{rem:smoothly trivial})
$$\Em = \{e : X \lto \C : e \text{ is smooth }\}.$$
It is well-known that $\cH^1(\Os_h) = \{0\}$ (\S 3, Chapter 0, \cite{GH1}).  Hence $\cM_h$ consists of a single point.

\subsubsection{The smooth de Rham moduli space}
A smooth integrable connection $\nb$ on $\Em$ is of the form $d + \eta$ where $\eta$ is a closed 1-form in $\Omega^1(\End(\Em)) \cong \Omega^1$.  The de Rham moduli space $\MdR = \ker(d)/\im(d\log)$ is the sheaf cohomology of
$$\G \stackrel{d\log}{\lto} \Omega^1 \stackrel{d}{\lto} \Omega^2,$$
where $\G$ is the gauge group consisting of maps $g : X \lto \C^\times$.  Hence it induces a map on the fundamental group $g_* : \pi_1(X) \lto \pi_1(\C^\times)$ which is an endomorphism of $\Z$.  In other words $g(z) = z^k g'(z)$, where $g_*' : \pi_1(X) \lto \pi_1(\C^\times)$ is the trivial map.  By homotopy lifting, we can consider $g' \in \Omega^0 = \Os.$  In other words,  $\G \cong \Z \times \Os.$  There is a decomposition
$$\Omega^1 = \{\frac{a dz}{z}\}\oplus d(\Os), \ \ a \in \C.$$
The de Rham cohomology $\H_d^1(\Os) \cong \C$ is generated by the cocycle $[\frac{dz}{z}]$ and $\MdR \cong \C/\Z$ is a torus as expected.  One can also obtain this from
the de Rham theorem which states that $\H_d^1(\Os) \cong \cH^1(\Os_c)$.

\subsubsection{The holomorphic structures} A holomorphic structure on $\Em$ is also defined by $\db + \eta$ where $\eta \in \Omega^{0,1}(\End(\Em)) $.  The space $\M_h$ is the sheaf cohomology of
$$\G \stackrel{\dbar\log}{\lto} \Omega^{0,1} \stackrel{\dbar}{\lto} \Omega^{0,2}.$$
Since $X$ has dimension 1, $\Omega^{0,2} = \{0\}.$
One can compute this directly or use the
Dolbeault theorem which states that $\H_\dbar^{0,1}(\Os) \cong \cH^1(\Os_h)$.  Thus $\M_h \cong \cH^1(\Os_h^\times)$ and consists of a single point.

\subsubsection{Holomorphic integrable connections ($\lambda = 1$)}
Since there is only one holomorphic structure up to isomorphism, we assume that it is $\dbar$.  The moduli space of integrable connections is the sheaf cohomology of
$$\G_h \stackrel{\d\log}{\lto} \Omega_h^1 \stackrel{\d}{\lto} \Omega_h^2.$$
Since $X$ is of complex dimension 1, $\Omega_h^2 = \{0\}$.
The holomorphic gauge group $\G_h$ consists of holomorphic maps $g : X \lto \C^\times$.  Again, the induced map $g_* : \pi_1(X) \lto \pi_1(\C^\times)$ is an endomorphism on $\Z$.  Hence $g(z) = z^k g'(z)$, where $g_*' : \pi_1(X) \lto \pi_1(\C^\times)$ is the trivial map.  By homotopy lifting, one may assume that $g' \in \Omega_h^0 \cong \Os_h.$  In other words,  $\G_h \cong \Z \times \Os_h$.

To further elaborate, the $\Os_h$-module of closed holomorphic 1-forms is of the form $f(z)dz$ where $f$ is a holomorphic function on $X$. Let $A = \frac{a dz}{z} \in \Omega_h^1$.  Since $G$ is abelian, $\Fsf(\nb) = \d A = 0$ and $\nabla = \d + A$ is an integrable connection.  However $A$ is not exact because the anti-derivative of $\frac{1}{z}$ is $\log(z)$ which is not defined on the entire $X$.
The equation
$\nabla f = (\d + A) f = 0$ has solutions of the form $f = C z^{-a}$ which is multi-valued on $X$.   Hence the equation has no global solution on $X$.  This is reminiscent of the path dependence of formula~(\ref{eq:path_function}), except the dependence is only on isotopy classes.
The multi-valued function $f$ on $X$ is then by definition a function $f : \tX \lto \C$.  The $\pi$-action on $\tX$ induces a $\pi$-action on $p^{-1}(x)$, hence, an action on $\{f(y) : y \in p^{-1}(x)\}$, hence, a (monodromy) representation $$\pi \lto \Aut(\C) \cong \C^\times.$$
Since $X$ is a punctured disk, $\pi \cong \Z$. Therefore
the representation variety $\Hom(\pi_1(X), \C^\times)$ is also $\C^\times$.  In the end, $\MdR^1 \cong \C/\Z$ just as before.
\subsubsection{The Dolbeault moduli space ($\lambda = 0$)}  The Dolbeault moduli space $\MdR^0$ is parameterized by the trivial holomorphic structure $\db$ and $\cH^0(\Omega_h^1)$.  The gauge group acts by the adjoint action, hence, trivially.  Thus $\MdR^0 \cong \cH^0(\Omega_h^1) \cong \C$ which is not diffeomorphic to $\MB \cong \C^\times$.

\subsection{Compact Riemann surfaces}
Now we assume that $X$ is a compact Riemann surface of genus $g$.

\subsubsection{The Betti moduli space}  This is already done in Section~\ref{sec:Hom} with $G = \C^\times$.  The commutator relation is trivial; hence, $\MB = (\C^\times)^{2g}$.

\subsubsection{The \v{C}ech moduli space}

By the long exact sequence (\ref{eq:constant long exact sequence}), the moduli space of rank-1 locally constant sheaves is $$\cM_c = \ker{\delta} \cong \cH^1(\Os_c)/\cH^1(\Z) \cong \C^{2g}/\Z^{2g} \cong (\C^\times)^{2g}.$$
Similarly,
$$\cM_h = \ker{\delta} \cong \C^g/\Z^{2g}$$ is the moduli space of rank-1 holomorphic sheaves that embed into the trivial smooth sheaf (See Remark~\ref{rem:smoothly trivial}).
The space $\cM_h$ is the Jacobi variety of $X$.  Observe that the complex structure on $\MC_h$ depends on how $\Z^{2g}$ injects into $\C^g$.  This in turn depends on the complex structure on $X$.  This is in direct contrast to $\cM_c$ where the embedding $\Z^{2g} \hookrightarrow \C^{2g}$ is canonical and independent of the complex structure on $X$.

\subsubsection{The smooth de Rham moduli space}
A smooth integrable connection $\nb$ on $\Em$ is of the form $d + \eta$ where $\eta \in \Omega^1(\End(\Em))(X) \cong \Omega^1(X)$.  The de Rham moduli space $\MdR = \ker(d)/\im(d\log)$ is the sheaf cohomology of
$$\G \stackrel{d\log}{\lto} \Omega^1 \stackrel{d}{\lto} \Omega^2,$$
where $\G$ is the gauge group consisting of maps $g : X \lto \C^\times$.
The de Rham theorem states that $\H_d^1(\Os) \cong \cH^1(\Os_c)$ and from this, we conclude that
$$\MdR \cong \C^{2g}/\Z^{2g} \cong (\C^\times)^{2g}.$$

\subsubsection{The holomorphic structures} A holomorphic structure on $\Em$ is also defined by $\db + \eta$ where $\eta \in \Omega^{0,1}(\End(\Em)) $.  The space $\M_h$ is the sheaf cohomology of
$$\G \stackrel{\dbar\log}{\lto} \Omega^{0,1} \stackrel{\dbar}{\lto} \Omega^{0,2}.$$
Since $X$ has complex dimension 1, $\Omega^{0,2} = \{0\}.$
The
Dolbeault theorem states that $\H_\dbar^{0,1}(\Os) \cong \cH^1(\Os_h)$.  The moduli space of holomorphic structure is
$$\M_h \cong \cM_h \cong \ker(\delta) \cong \C^g/\Z^{2g}.$$
Again this is the Jacobi variety of $X$.

\subsubsection{Holomorphic integrable connections ($\lambda = 1$)}
An object in this groupoid consists of a pair $(\bn, \nb^1)$ where $\bn \in \Fh$ and $\nb^1$ is a holomorphic integrable connection on $\Em_h \cong \Em^\bn$.

However, with this construction, there is a canonical rational map $\MdR^1 \dashrightarrow \M_h$ arising from (See (\ref{eq:projection to the graded}))
$$\FF_d^1 \lto \Fh, \ \ (\bn,\nb^1) \mapsto \bn.$$
In this case of $G = \C^\times$, this rational map is actually defined every where, hence, holomorphic.  To be more explicit, there is a canonical holomorphic map
$$
\MdR^1 \cong (\C^\times)^{2g} \lto \C^g/\Z^{2g} \cong \M_h \cong \cM_h.
$$
Recall again that $\M_h$ is actually the Jacobi variety whose complex structure depends on the complex structure on $X$ while the complex structure on $\MdR^1$ does not depend on the complex structure on $X$.

\subsubsection{The Dolbeault moduli space ($\lambda = 0$)}  The Dolbeault moduli space $\MdR^0$ parameterizes pairs $(\bn, \nb^0)$ where $\bn \in \Fh$ and $\nb^0 \in \Omega_h^1(X)$. Hence there are respectively canonical holomorphic and rational maps $\MdR^0 \stackrel{\dashrightarrow}{\leftarrow} \M_h$ arising from
$$\Fh \lto \FF_d^0, \ \ \bn \mapsto (\bn, 0) \ \text{ and } \ \FF_d^0 \lto \Fh, \ \ (\bn,\nb^0) \mapsto \bn.$$
The general theory states that when $X$ is projective, there are the following topological (diffeomorphic) equivalences
$$\MdR^0 \ \cong \MdR^1 \cong \MdR \cong \MC_c.$$
In this case of $G = \C^\times$, the general theory effectively reduces to the classical Hodge theory.  Recall that the classical Hodge theorem states that
$$
\H_d^1(\Os) \cong \H_\dbar^{0,1}(\Os) \oplus \H_\dbar^{1,0}(\Os) \cong \cH^1(\Os_h) \oplus \cH^0(\Omega_h^1).
$$
From this, we obtain a direct sum decomposition
$$\MdR^0 \cong \M_h \oplus \H_\dbar^{1,0}(\Os) \cong \cM_h \oplus \cH^0(\Omega_h^1).$$
We want to emphasize here that the identification of $\MdR^0$ and $\MdR$ is homeomorphic (diffeomorphic), but neither algebraic nor holomorphic -- $\MdR^0$ contains a compact complex torus while $\MdR$ does not.
\section{Some final remarks}\label{sec:final remarks}
The moduli spaces $\MB$ and  $\MdR$ have the same underlying topology.  When $X$ is smooth projective, so does the Dolbeault moduli space $\MdR^0$.  Denote by $\M$ the underlying topological space.
When $G$ is non-abelian, the topology of $\M$ is complicated.  One can see this already from the Betti moduli space with its inherited variety structure from $G$.  When $X$ is smooth projective, the Dolbeault moduli space $\MdR^0$ is homeomorphic (diffeomorphic) to $\MdR$, hence, to $\MB$ \cite{Co1, Do1, GM1, GX1, Hi1, Si0}.  This gives an extra complex structure on $\M$.  The smooth part of the Betti moduli space $\MB$ also has a natural underlying symplectic structure \cite{G1} which, together with the two complex structures on $\M$, give rise to a hyperk\"ahler structure on $\M$.  This hyperk\"ahler structure yields rich topological information of $\M$ \cite{Hi1} and the most pleasant surprise occurs when the structure group is real:

Let $K = \U(r) \subseteq \GL(r,\C) = G$, the maximum compact subgroup of $G$.  Readers can immediately check that the constructions for the Betti, the smooth and constant \v{C}ech and the smooth de Rham groupoids are still valid.  Denote them as $\Hom(\pi,K), \cF(K), \cF_c(K), \FF_d(K)$, respectively.
The inclusion $P : \cF_c \lto \cFh$ of Remark~\ref{rem:underlying holomorphic structure} on objects is injective, but the induced map on the isomorphism classes $P : \Iso(\cF_c) \lto \Iso(\cFh)$ may neither be injective nor surjective.

Continuing from last chapter, we now assume that $X$ is a compact Riemann surface of genus $g$ and $r=1$, i.e. $K = \U(1) \subseteq \C^\times = G$.  Then $\cM_h \cong \C^g/\Z^{2g}$ is a complex torus and so is $\MB(K) = \Hom(\pi,K)/K$ and the composition map
$$ \Iso(\cF_c(K)) \lto \Iso(\cF_c) \stackrel{P}{\lto} \Iso(\cFh)$$
is one-to-one and onto.
More specifically, this gives the following diffeomorphisms
$$\U(1)^{2g} \cong \Hom(\pi,\U(1))/\U(1) \cong \cM_h.$$
This provides $\MB(\U(1)) :=  \Hom(\pi,\U(1))/\U(1)$ with a complex structure.

In general, for a smooth complex projective variety $X$ and a compact Lie group $K$, the representation variety $\MB(K)$ acquires the complex variety structure of $\M_h$ \cite{Hi1, Si3}.  It so happens that the same is true for any real form $H \subset G$.  That is if $H$ is a real form of $G$, then $\MB(H) := \Hom(\pi_1(X), H)^+/H$ acquires a complex variety structure \cite{Hi1, Si3}.
Notice that in the case of $G = \C^\times$, there are diffeomorphisms
$$\MB(\U(1)) \oplus \MB(\R^\times) \cong \MB \cong \cM_h \oplus \cH^0(\Omega_h^1).$$
To summarize, if $G$ is complex reductive, then $\M$ has a hyperk\"ahler structure; moreover, if $H \subset G$ is a real form, then $\M(H)$ (with underlying topological space $\MB(H)$) acquires a complex variety structure.
\vskip 0.1in
\noindent {\bf Acknowledgment:} The author thanks Nan-Kuo Ho for proof-reading the manuscript.


\begin{thebibliography}{99}
\bibitem{Bo1} Borel, A.; Grivel, P.-P.; Kaup, B.; Haefliger, A.; Malgrange, B.; Ehlers, F. Algebraic $D$-modules. Perspectives in Mathematics, 2. {\em Academic Press, Inc., Boston, MA}, 1987.

\bibitem{Co1} Corlette, Kevin Flat $G$-bundles with canonical metrics.  {\em J. Differential Geom.}  {\bf 28}  (1988),  no. 3, 361--382.

\bibitem{Co2} Coutinho, S. C. A primer of algebraic $D$-modules. London Mathematical Society Student Texts, 33. {\em Cambridge University Press, Cambridge}, 1995. xii+207 pp.

\bibitem{Do1} Donaldson, S. K.
Twisted harmonic maps and the self-duality equations.
{\em Proc. London Math. Soc. (3)} {\bf 55} (1987), no. 1, 127--131.

\bibitem{G0} Goldman, William M. Higgs bundles and geometric structures on surfaces. {\em The many facets of geometry,} 129–-163, {\em Oxford Univ. Press, Oxford,} 2010.

\bibitem{G1} Goldman, William M. The symplectic nature of fundamental groups of surfaces. {\em Adv. in Math.} {\bf 54} (1984), no. 2, 200--225.

\bibitem{GM1} Goldman, William M.; Millson, John J. The deformation theory of representations of fundamental groups of compact K\"{a}hler manifolds.  {\em Inst. Hautes \'{E}tudes Sci. Publ. Math. No.67}  (1988), 43--96.

\bibitem{GX1} Goldman, William M.; Xia, Eugene Z. Rank one Higgs bundles and representations of fundamental groups of Riemann surfaces.  {\em Mem. Amer. Math. Soc.}  {\bf 193}  (2008),  no. 904, viii+69 pp.

\bibitem{GH1} Griffiths, Phillip; Harris, Joseph Principles of algebraic geometry. Pure and Applied Mathematics. {\em Wiley-Interscience [John Wiley \& Sons], New York}, 1978. xii+813 pp.

\bibitem{Ha1} Hartshorne, Robin Algebraic geometry. Graduate Texts in Mathematics, No. 52. {\em Springer-Verlag, New York-Heidelberg}, 1977. xvi+496 pp.

\bibitem{Hi1} Hitchin, N. J. The self-duality equations on a Riemann surface.  {\em Proc. London Math. Soc. (3)}  {\bf 55}  (1987),  no. 1, 59--126.

\bibitem{Ka1} Katz, Nicholas M. Rigid local systems. Annals of Mathematics Studies, {\bf 139}. Princeton University Press, {\em Princeton, NJ,} 1996. viii+223 pp.

\bibitem{Mu1} Mumford, D.; Fogarty, J.; Kirwan, F. Geometric invariant theory. Third edition. Ergebnisse der Mathematik und ihrer Grenzgebiete (2) [Results in Mathematics and Related Areas (2)], {\bf 34}. {\em Springer-Verlag, Berlin}, 1994. xiv+292 pp.

\bibitem{NN1} Newlander, A.; Nirenberg, L.
Complex analytic coordinates in almost complex manifolds.
Ann. of Math. (2) {\bf 65} (1957), 391--404.

\bibitem{Ne1} Newstead, P. E. Introduction to moduli problems and orbit spaces. Tata Institute of Fundamental Research Lectures on Mathematics and Physics, {\bf 51}. {\em Tata Institute of Fundamental Research, Bombay; by the Narosa Publishing House, New Delhi}, 1978. vi+183 pp.

\bibitem{Ri1} Riemann, Georg Friedrich Bernhard Deux the\'eor\`emes g\'en\'eraux sur les \'equations diff\'erentielles lin\'eares a coefficients alg\'ebriques, Œuvres mathématiques de Riemann. (French) Traduites de l'allemand par L. Laugel, avec une préface de C. Hermite et un discours de Félix Klein. Nouveau tirage {\em Librairie Scientifique et Technique Albert Blanchard, Paris} 1968, 353-368.

\bibitem{Sc1} Schneiders, Jean-Pierre An introduction to $D$-modules. Algebraic Analysis Meeting (Lie\`ge, 1993). {\em Bull. Soc. Roy. Sci. Li\`ege} {\bf 63} (1994), no. 3-4, 223--295

\bibitem{Si0} Simpson, Carlos T. Constructing variations of Hodge structure using Yang-Mills theory and applications to uniformization. {\em J. Amer. Math. Soc.} {\bf 1} (1988), no. 4, 867--918.

\bibitem{Si1} Simpson, Carlos T.
Higgs bundles and local systems.
{\em Inst. Hautes \'Etudes Sci. Publ. Math. No. 75} (1992), 5--95.

\bibitem{Si2} Simpson, Carlos T. Moduli of representations of the fundamental group of a smooth projective variety. I.  {\em Inst. Hautes \'{E}tudes Sci. Publ. Math. No. 79}  (1994),   47--129.

\bibitem{Si3} Simpson, Carlos T. Moduli of representations of the fundamental group of a smooth projective variety. II.  {\em Inst. Hautes \'{E}tudes Sci. Publ. Math. No. 80}  (1995),  5--79.

\end{thebibliography}
\end{document}